\documentclass{amsart}
\usepackage{amsthm,amstext,amscd,amsbsy,amsxtra,latexsym}

\usepackage{cite}

\usepackage[citation-order, msc-links]{amsrefs}

\usepackage{graphicx}
\usepackage{url}
\usepackage{hyperref}

\hypersetup{
  colorlinks   = true, %Colours links instead of ugly boxes
  urlcolor     = blue, %Colour for external hyperlinks
  linkcolor    = blue, %Colour of internal links
  citecolor   = red %Colour of citations
}

 \usepackage[usenames,dvipsnames]{color}

\usepackage[capitalise]{cleveref}
\crefname{thm}{Thm.}{}
\crefname{prop}{Prop.}{}
\crefname{lem}{Lem.}{}
\crefname{cor}{Cor.}{}
\crefname{prob}{Problem}{}
\crefname{figure}{Fig.}{}

\usepackage{color}

 \newcommand{\ts}{{{\tilde{s}}}}
 \newcommand{\I}{{\rm{i}}}

 \newcommand{\Z}{{\mathbb Z}}

 \newcommand{\Q}{{\mathbb Q}}
 \newcommand{\Pp}{{\mathbb P}}
 \newcommand{\CC}{{\mathbb C}}

 \newcommand{\Ee}{{\mathcal E}}

 \newcommand{\KK}{{\mathcal K}}

 \newcommand{\Sc}{{\mathcal S}}

 \newcommand{\gal}{{\text{Gal}}}

 \newtheorem{thm}{Theorem}[section]

 \numberwithin{equation}{section}

\newcommand\iso{\cong}

\def\det{\mbox{det }}

\def\deg{\mbox{deg }}

\usepackage{color}
\usepackage{booktabs}

\begin{document}
	
	\title{Generators and splitting fields of certain elliptic K3 surfaces}
	\author{Sajad Salami}
	\address{Institute of Mathematics and Statistics, State University of Rio de Janeiro, Rio Janeiro, Brazil}
	\email{sajad.salami@ime.uerj.br}
	
	\author{Arman Shamsi Zargar}
	\address{Department of Mathematics and Applications, University of Mohaghegh Ardabili, Ardabil, Iran}
	\email{zargar@uma.ac.ir}
	
	\maketitle

	\begin{abstract}
	Let $k \subset {\mathbb C}$ be a number field and ${\mathcal E}$ be an elliptic curve  defined over $k(t)$,  the  rational  function field  of the projective line
	${\mathbb P}^1_k$,
	is  isomorphic to the generic fiber of an elliptic surface  
	$\pi:= \Sc_\Ee \rightarrow  {\mathbb P}^1_k$.  
	For any subfield  ${\mathcal K}\subseteq {\mathbb C}$ of $k$, the set ${\mathcal E}({\mathcal K}(t))$ of ${\mathcal K}(t)$-rational points of ${\mathcal E}$ is known to be a finitely generated abelian group.
	The  splitting field of ${\mathcal E}$ defined over $k(t)$ is the smallest finite extension ${\mathcal K} \subset {\mathbb C}$ of $k$  such that  
	${\mathcal E} ({\mathbb C} (t)) \iso {\mathcal E} ({\mathcal K}(t))$. 
	In this paper, we consider the elliptic $K3$ surfaces defined over  $k={\mathbb Q}$ with the generic fiber given by  the Weierstrass  equation  ${\mathcal E}_n:  \displaystyle y^2=x^3 + t^n +  1/t^n$, $1\leq n\leq 6$, and determine the splitting field ${\mathcal K}_n$, and find an explicit set of independent generators for   ${\mathcal E}_n ({\mathcal K_n}(t))$  for  $1\leq n \leq 6$.
		\end{abstract}
	
%	\begin{keyword}
%		Elliptic surface, Mordell-Weil lattice, Splitting  field.
%		\MSC[2020]  14J27  \sep    11G05. 
%	\end{keyword}
	
	\maketitle
%	\tableofcontents

\section{Introduction and main results}

Let $k$ be a number field and 
${\mathcal E}$ be an elliptic curve defined over   $k(t)$,  the rational function field of
the  projective line $\Pp^1_k$  over $k$,  
that is isomorphic to the generic fiber of an elliptic surface   
$\pi:= \Sc_\Ee \rightarrow  {\mathbb P}^1_k$.
Given any subfield $\KK \subseteq \CC(t)$,  the set $\Ee(\KK)$ of   $\KK$-rational points of ${\mathcal E}$ is known to be a finitely generated abelian group and has a lattice structure called the Mordell--Weil lattices  \cite{Shioda1990a, Silverman1994, Schuett2019}.

By 	the \textit{splitting field}
of $\Ee$ over $k(t)$, we mean   the smallest finite extension $\KK\subset \CC$ of $k$, 
such that  $\Ee (\CC (t))\iso \Ee (\KK(t))$.
It is a well-known fact  that  $\KK| k$ is a Galois  extension with the finite Galois group  $G=\gal(\KK | k)$. 
Moreover, 	  the $G$-invariant elements of $\Ee (\KK(t))$ are the $\Ee (k(t))$-rational points \cite{Schuett2019}.

%In \cite{Shioda1990a}, T.~Shioda introduced the theory of Mordell--Weil  lattices associated to an elliptic curve $\Ee$ defined over $\CC(t)$,  which is known to be  isomorphic	with $\Ee (\CC (t))$.	 The  splitting field can also be defined in relation with  Mordell--Weil lattices as follows. There is a natural action of $\gal(\CC | k)$ on $\Ee(\CC (t))$ that preserves the height pairing $\left\langle ,  \right\rangle$ and gives a Galois representation
%$$\rho :  \gal(\CC | k) \longrightarrow  \aut (\Ee (\CC (t)), \left\langle ,  \right\rangle ).$$
%Since  the height pairing   $\left\langle ,  \right\rangle  $ on $\Ee (\CC(t))$
%is positive definite up to torsion subgroups, the image ${\rm Im}(\rho)$  of $\rho$ is a subgroup of  the finite group $\aut (\Ee (\CC (t)), \left\langle ,  \right\rangle ).$
%Hence, in the terminology of  Galois theory, the splitting field $\KK$ is exactly  the extension
%of $k$ which corresponds to  $\ker (\rho)$, and  we have  $\gal(\KK | k)= {\rm Im}(\rho).$
%One can see more on the general theory of Mordell--Weil lattices in  \cite{Schuett2019, Shioda1990a, Silverman1994}
%and  associated Galois representations in \cite{Schuett2019,  Shioda1989c, Shioda1989b, Shioda1989a}.

In this paper, we  consider $k=\Q$ and  the elliptic $K3$ surfaces 	 over $\Q(t)$ with  a generic fiber given by the  following equation 
\begin{equation*}
	%\label{eq0}
	\Ee_n : y^2=x^3 + t^n + \frac{1}{t^n}, \ \text{for} \ 1 \leq n \leq 6.
\end{equation*}
The structure of  Mordell-weil lattice  of $\Ee_n$ over $\CC(t)$ is studied by T.~Shioda in \cite{Shioda2000a, Shioda2008} and by A.~Kumar and M.~Kuwata in \cite{Kumar2017} with a more general setting, for all $1 \leq n \leq 6$.
We notice  that  ${\mathcal E}_n$ is a special member, considering  $\alpha=\beta=0$,
of  the  generic fiber of a more general family $K3$ surface   
defined  by 
$$%\Ff^{(n)}_{\alpha,\beta}: 
y^2=x^3 -3 \alpha x + \left( t^n + \frac{1}{t^n} - 2 \beta \right).$$
In particular, the invariants of the Mordell--Weil lattices of $\Ee_n$ are determined  by T.~Shioda in    \cite[Theorem~2.4]{Shioda2008}, and
a generic form of  their   generators is described in    \cite[Theorem~2.6]{Shioda2008}.
For the convenience of readers, we gathered  those results  as   Theorem~\ref{shi-thm} in Section~\ref{prel}.

The main aims of this paper are to determine  the splitting field $\KK_n \subset \CC$  of $\Ee_n$ and  provide  an explicit set of  independent  generators of  $ \Ee_n(\KK_n(t))$ for each  $1\leq n \leq 6$.  Let  $r_n$  be the rank of
$ \Ee_n(\CC(t))$, and  $\zeta_m$  be  a fixed $m$-roots of unity.
In table  \ref{tab:notation}, we gathered all  $m$-roots of unity  and  algebraic quantities
that  we used throughout the paper.

\begin{table}[htbp]
	\caption{Notations}\label{tab:notation}
	\begin{tabular}{|l|l|l|}
		\hline
		%		$n$ &  $1$ & $2$ \\ 
		%		\midrule
		$\I= \sqrt{-1}$ & $\epsilon_1=2+ \sqrt{3}$ &    $\epsilon'_1=2- \sqrt{3}$  \\ [5pt]
		$\displaystyle\zeta_3 = \frac{\mathrm{i} \sqrt{3} -1}{2}$ & $\epsilon_2  = 11 \sqrt{2} + 9 \sqrt{3}$  &    $\epsilon'_2  = -11 \sqrt{2} + 9 \sqrt{3}$   \\ [5pt]
		$\displaystyle\zeta_5 =\frac{\sqrt{5}-1 + \I \sqrt{2} \sqrt{5 +\sqrt{5}}}{4}$ &  $\epsilon_3  =  \sqrt{2} + 5 \I$ & $\epsilon'_3  =  \sqrt{2} - 5 \I$  \\ [5pt]
		$\displaystyle\zeta_6 =  \frac{1+\mathrm{i} \sqrt {3} }{2}$ &  $\epsilon_4 =1-\zeta_{12}$	& $\beta_0= 2^{\frac{1}{6}}$ \\ [5pt]
		$\displaystyle\zeta_8 =  \frac{\sqrt{2} ( 1+\mathrm{i} )}{2}$ &  $\displaystyle \epsilon_5 =\frac{\sqrt {3}-\sqrt {5}}{2}(\zeta_{12}+ \zeta_{12}^{10})$ &
		$\beta_1= \left(3+2 \sqrt{3}\right)^{\frac{1}{4}}$   \\ [5pt] 
		$\displaystyle\zeta'_8 =\frac{\sqrt{2} ( 1-\mathrm{i} )}{2}$ & $\displaystyle \epsilon_6  =\frac{1 + \sqrt {5}}{2} \zeta_{12}$  & $\beta_2=\left(3-2 \sqrt{3}\right)^{\frac{1}{4}}$ \\ [5pt]
		$\displaystyle\zeta_{12} = \frac{\I + \sqrt {3}}{2}$ &  & \\ [5pt]
		\hline
	\end{tabular}
\end{table}

The Mordell--Weil lattice   $\Ee_1(\CC(t))$ is of rank $r_1=0$, see \cite[Theorems~1.1]{Shioda2007b}.
In the rest of this section, we	provide a list of our main results for each of the cases $2\leq n \leq 6$. 

\begin{thm}
	\label{main-0}
	The Mordell--Weil lattice   $\Ee_2(\CC(t))$   is isomorphic to $\Ee_2(\KK_2(t)) $  with $r_2=4$, where
	$\KK_2= \Q(\zeta_3, 2^{\frac{1}{3}})$  of degree 6 with a minimal  defining polynomial $ g_2(x)=x^6+ 108.$
	
	Moreover,  a set of   linearly independent generators
	of $\Ee(\KK_2(t))$ includes the following four points:
	$$\begin{aligned}[b]
		P_1 &=\left(2^{\frac{1}{3}},  t+\frac{1}{t}\right),  & 	P_2 &=\left(  2^{\frac{1}{3}}\zeta_3,  t+\frac{1}{t} \right), \\ 
		P_3 & =\left( -2^{\frac{1}{3}}, t-\frac{1}{t}\right),  & P_4 & =\left(- 2^{\frac{1}{3}}\zeta_3^2, t-\frac{1}{t} \right).
	\end{aligned}$$
\end{thm}

\begin{thm}
	\label{main-a}		
	The Mordell--Weil lattice   $\Ee_3(\CC(t))$   is isomorphic to $\Ee_3(\KK_3(t)) $   with $r_3=8$, 
	where $\KK_3 = \Q(\zeta_3, (3+ 3 \sqrt{3})^{\frac{1}{4}})$
	with a minimal  defining $g_3(x)$ of degree 16 given by \ref{f3}.
	% and contains 
	%	$\Q \left(\zeta_{12},   3^{\frac{1}{8}}, 
	%	\epsilon_1^{\frac{1}{4}} \right).$
	
	Moreover, a set  of eight independent generators of $\Ee(\KK_3(t))$ includes the following  points:
	%		$P_j=(x_j(t), y_j(t) )$ and $P_{j+8}=(x_{j+8}(t), y_{j+8}(t))$,
	$$	P_j = \left( x_j(t), y_j(t)\right) = \left( \frac{a_j t^{2}+b_j t +a_j}{t}, \ \frac{c_j t^{2}+d_j t +c_j}{t} \right),$$
	and $P_{j+4}=\zeta_3 \cdot P_{j}= \left( x_{j}(\zeta_3 t), y_{j}(\zeta_3 t)\right) $    for $j=1, 2, 3, 4$,
	%	\begin{align*}
		%		P_j   & = (x_j(t), y_j(t)= \left( \frac{a_j \,t^{2}+b_j t +a_j}{t}, \ \frac{c_j \,t^{2}+d_j t +c_j}{t} \right),\\
		%		P_{j+4} & = \left(\frac{a_j \zeta_3^{2} \,t^{2} +b_j \zeta_3 t +a_j}{\zeta_3 t}, \
		%		\frac{c_j \zeta_3^{2}\,t^{2}+d_j \zeta_3 t +c_j}{\zeta_3 t} \right),
		%	\end{align*}
	where   $a_j, b_j, c_j, d_j$ are given in Subsection~\ref{coef3}.
\end{thm}
For $n= 4, 6$, we consider the automorphism   of $\Ee_n(\CC(t))$  given by
$$\phi_n: (x(t), y(t))\rightarrow (-x( \zeta_{2n}  t), \I  y (\zeta_{2n} t)).$$ 

\begin{thm}
	\label{main-b}
	The Mordell--Weil lattice   $\Ee_4(\CC(t))$ is isomorphic to $\Ee_4(\KK_4(t)) $   with $r_4=12$, where  $\KK_4 $ is   defined  by polynomial $g_4(x)$ of degree 24  given by \ref{g14}, containing  the number field
	$\Q \left(\zeta_{8}, \zeta_{12}, 2^{\frac{1}{12}}, \epsilon_2^{\frac{1}{6}}, \epsilon_3^{\frac{1}{6}} \right).$
	
	Moreover,  a set of  12 linearly  independent generators of $\Ee_4(\KK_4(t)) $ includes the following  points,
	$$P_j=\left( x_j(t), y_j(t)\right) =\left(\frac{a_j t^{2}+b_j t +a_j}{t}, 
	\frac{t^{4}+ c_j t^{3}+ (d_j +2) t^{2}+c_j t  +1}{t^{2}}\right),$$
	%	\begin{align*}
		%		 P_j&=\left( \frac{a_j \,t^{2}+b_j t +a_j}{t}, \, \frac{t^{4}+ c_j \,t^{3}+d_j \,t^{2}+c_j t +2 t^{2}+1}{t^{2}}\right)\\
		%		 P_{j+6} &= \left(- \frac{\left(1+\mathrm{i}\right) a_j \,t^{2}+2 b_j t +\left(2-2 \,\mathrm{i}\right) a_j}{2 t}, \, \frac{ \mathrm{i} t^{4}+\left(1+\mathrm{i}\right) c_j \,t^{3}+\left(4+2 d_j \right) t^{2}+\left(2-2 \,\mathrm{i}\right) c_j t -4 \,\mathrm{i}}{2 t^{2}}\right)
		%	\end{align*} 
	and  $P_{j+6}= \phi_4 (P_j)=\left(- x_j(\zeta_8 \, t), \I \, y_j(\zeta_8 \,  t)\right)$ for $j=1, \ldots, 6$,
	where    $a_j, b_j, c_j, d_j$ are given in Subsection~\ref{coef4}.
\end{thm}

\begin{thm}
	\label{main-c}
	The Mordell--Weil lattice   $\Ee_5(\CC(t))$ is isomorphic to $\Ee_5(\KK_5(t)) $   with  $r_5=16$, where $\KK_5=\KK'_5(\zeta_5)$and  $\KK'_5$ is a number field  of degree 96,	
	and    $\KK_5=\KK'_5(\zeta_5)$ has degree 192,  with   minimal  defining polynomials given in \cite[min-pols]{k3-codes}. 
	The splitting field  $\KK_5$ contains 
	the number field $\Q \left(\zeta_5,  \zeta_{12},  5^{\frac{1}{24}},  (\epsilon_4 \epsilon_5)^{\frac{1}{2}}\right),$
	where $\epsilon_4$ and $\epsilon_5$ 
	with
	$$	\begin{aligned}[b]\epsilon_4&=1-\zeta_{12}, & 
		\epsilon_5&=\left( \frac{1 + \sqrt {5}}{2}\right)  \zeta_{12}, &	
		\epsilon_5&=\left( \frac{\sqrt {3}-\sqrt {5}}{2}\right) (\zeta_{12}+ \zeta_{12}^{10}),
	\end{aligned}$$
	are the fundamental units of the number field $\Q(\I,\sqrt {3},\sqrt {5})=\Q(\zeta_{12},\sqrt{5})$.
	
	Moreover,  a set of sixteen independent generators of $\Ee_5(\KK_5(t)) $ includes $P_j=\left( x_j(t), y_j(t)\right) $ with
	$$\begin{aligned}
		x_j(t) &= \frac{t^4 + a_j t^3 + (b_j +2) t^2 + a_j t +1}{u_j^2  t^2}, \\
		y_j(t) &=\frac{t^6 +c_j t^5+ (d_j +3) t^4 +(2 c_j +e_j) t^3 +(d_j +3) t^2 + c_j t +1}{u_j^3 t^3},\\
	\end{aligned}
	$$
	%	and
	%	\begin{align*}
		%	x_{j+8}(t) &= \frac{\zeta_5^4\, t^4 + a_j\, \zeta_5^3\, t^3 + (b_j +2) \zeta_5^2\, t^2 + a_j \zeta_5\, t +1}{u_j^2 \, \zeta_5^2\, t^2}, \\
		%		y_{j+8}(t) &=\frac{\zeta_5\, t^6 +c_j t^5+ (d_j +3)\, \zeta_5^4\, t^4 +(2 c_j +e_j)\, \zeta_5^3\, t^3 +(d_j +3)\,\zeta_5^2\, t^2 + c_j\, \zeta_5\, t +1}{u_j^3 \,\zeta_5^3\, t^3},
		%	\end{align*}
	and $P_{j+8}=\left( x_{j}(\zeta_5 t), y_{j}( \zeta_5 t)\right) $  	for $j=1, \ldots, 8$,
	where   $a_j, b_j, c_j, d_j, e_j$  and  $u_j$'s are given in \cite[Points-5]{k3-codes}.
\end{thm}

We have to mentioned that the  points given by Theorem~\ref{main-c} provided sixteen   generators of the Mordell-Weil lattice of  the Shioda's rank 68 elliptic surface, as described in \cite[Thm. 1.1 (11)]{Salami-2025}.

\begin{thm}
	\label{main-d}
	The Mordell--Weil lattice   $\Ee_6(\CC(t))\cong \Ee_6(\KK_6(t))$ is isomorphic to $\Ee_6(\KK_6(t)) $   with $r_6=16$, where
	$\KK_6 $ is a number field with a defining minimal  polynomial $g_6(x)$ of degree 96
	%=\Q(\zeta_{8}, \zeta_{12},  2^{\frac{1}{12}},  3^{\frac{1}{8}}, \epsilon_1^{\frac{1}{4}},   {v_7}^{\frac{1}{2}}),$$ 
	%	 $$ \text{with} \	\beta_0= 2^{\frac{1}{6}},\,  \beta_1= (3+ 2 \sqrt{3})^{\frac{1}{4}}, $$
	%	and $v_7$ is given by \eqref{6eqn2}.
	given in \cite[min-pols]{k3-codes}.	
	Moreover,  a set of  16 independent generators includes $P_j=(x_j(t), y_j(t) ) $ with
	\begin{align*}
		x_j(t)& = \frac{a_{j,0} t^{4}+a_{j,1} t^{3}+a_{j,2} t^{2}+ a_{j,1} t +a_{j,0}}{t^2}, \\
		y_j(t) &=\frac{b_{j,0} t^{6}+b_{j,1} t^{5}+b_{j,2} t^{4}+b_{j,3} t^{3}+b_{j,2} t^{2}+b_{j,1} t +b_{j,0}}{t^3},
		%\end{align*}
		%	
		%	\begin{align*}
			%a_{j,\,4}& =a_{j,\,0}=a_j, \ \ a_{j,\,3} = a_{j,\,1} =b_j-2 \sqrt{2}\,  a_j,  \ \
			%a_{j,\,2} =  g_j+4 a_j -\sqrt{2}\, b_j,
			%\\
			%b_{j,\,6}&=b_{j,\,6}=c_j, \ b_{j,\,5}=b_{j,\,1}= c_j +d_j -3 \sqrt{2}, \\
			%b_{j,\,4}&=b_{j,\,2}= d_j +9 c_j +e_j -2 \sqrt{2}, \
			%b_{j,\,3}=h_j-8 \sqrt{2}\, c_j -\sqrt{2}\, e_j +4 d_j,
		\end{align*}
		in which
		$$\begin{aligned}[b]
			a_{j,0}& =a_j, & a_{j,1} &=b_j-2 \sqrt{2}  a_j,  \\
			a_{j,2} &=  g_j+4 a_j -\sqrt{2} b_j, &  &
			\\
			b_{j,0}&=c_j, & b_{j,1}&= c_j +d_j -3 \sqrt{2}, \\
			b_{j,2}&= d_j +9 c_j +e_j -2 \sqrt{2}, &
			b_{j,3}&=h_j-8 \sqrt{2} c_j -\sqrt{2} e_j +4 d_j, &&
		\end{aligned}$$
		and the points $P_{j+8}= \phi_6(P_j)$  	for $j=1, \ldots, 8$,
		%	\begin{align*}
			%	x_{j+8}(t)& = \frac{a_{{j+8},\,4} \,t^{4}+a_{{j+8},\,3} \,t^{3}+a_{{j+8},\,2}\, t^{2}+ a_{{j+8},\,1}\, t +a_{{j+8},\,0}}{\zeta_{12}^2 \, t^2}, \\
			%		y_{j+8}(t) &=\frac{b_{{j+8},\,6}\,t^{6}+b_{{j+8},\,5}\, t^{5}+b_{{j+8},\,4}\, t^{4}+b_{{j+8},\,3}\, t^{3}+b_{{j+8},\,2}\, t^{2}+b_{{j+8},\,1}\, t +b_{{j+8},\,0}}{\zeta_{12}^3 \,  t^3},
			%	\end{align*}
		%where
		%		$$\begin{aligned}[b]
			%		a_{{j+8},\,4} & =\zeta_3 a_{{j+8},\,0} = \zeta_3\,  a_j,&
			%		a_{{j+8},\,3} & =\zeta_6 a_{{j+8},\,1} = \zeta_6\, 2 \sqrt{2}\, a_j+b_j,\\
			%		a_{{j+8},\,2} & = \zeta_6\, ( a_j + \sqrt{2}\, b_j + g_j),&
			%		b_{{j+8},\,6} & = -b_{{j+8},\,0} = c_j,\\
			%		b_{{j+8},\,5} & =b_{{j+8},\,1}=\zeta_{12}^5 (3 \sqrt{2} c_j + d_j),&
			%		b_{{j+8},\,4} & =b_{{j+8},\,2} =  \zeta_3 ( 9 c_j +2 \sqrt{2} d_j + e_j ),\\
			%		b_{{j+8},\,3} & = \mathrm{i} (8\sqrt{2} c_j + 4 d_j + \sqrt{2} e_j + h_j ), &
			%	\end{aligned}$$
		where   $a_j, b_j, c_j, d_j, e_j, g_j$ and $h_j$
		are given in \cite[Points-6]{k3-codes}.
	\end{thm}
	
	We would like to mention that the results of this paper are used in our under progress works. In  \cite{S-Shz-2023},  we attempt to explicitly determine the generators and splitting fields of the Shioda elliptic surface given by $y^2=x^3+ t ^m +1$ for integers $2 \leq m\leq 12$ defined over $\Q(t)$;
	and in  \cite{Salami-2025} for the particular case $m=360$, which is known to have rank 68 over ${\mathbb C}(t)$.
	
	In our computations, we mostly used 
	%the packages {\sf PolynomialTools} and {\sf PolynomialIdeals} in
	the mathematical software  {\sf Maple} \cite{maple},  and 
	{\sf Pari/Gp} \cite{PARI2}
	% as well as  {\sf{SageMath}}~\cite{sagemath}.
	
	The rest of  paper is  organized as follows.
	Prior to proving the main results, we provide the preliminary facts on the Mordell--Weil lattice of
	${\mathcal E}_n$ in the next section.
	Then, we prove   Theorems~\ref{main-0} and \ref{main-a} in Section~\ref{case123}.
	In the last three  sections, we respectively demonstrate the proof of Theorems~\ref{main-b}, \ref{main-c} and  \ref{main-d}.

	%==========================================
	\section{Shioda's results on Mordell--Weil lattice of  $\Ee_n$}
	\label{prel}
	In this section, we   recall some of the   known results by T.~Shioda on the elliptic $K3$ surfaces $\Ee_n$ defined over $\Q(t)$ from \cite{Shioda2008}. 
	
	For a given  lattice $(L, \left\langle ,  \right\rangle)$ and  an integer  $m\geq 2$,
	we let  $L[m]$  be a lattice with the height pairing $m \cdot \left\langle , \right\rangle.$
	We denote by $M_n$ the Mordell--Weil lattice $\Ee_n(\CC(t)),$
	which does not have  torsion part, see \cite[Lemma~5.2]{Shioda2008}.
	It is clear that $\Ee_n$ is obtained from $\Ee_1$ by the base change
	$t \rightarrow t^n$. Hence,  we let $N_n=M_1 [n]$  for  each $2\leq n \leq 6$.
	
	In order to  study the lattice  $M_n$, as in \cite{Shioda2008},  we will consider  Mordell--Weil lattice
	$ M'_n = \Ee'_n(\CC(s))$
	of the rational elliptic surface $\Ee'_n: y^2=x^3+ f_n(s)$
	% where $s=t+ 1/t$
	and
	$f_n(s)$ is a polynomial defined as follows,
	\begin{equation}
		\label{fns}
		f_n(s)=
		\begin{cases}
			s^2 -2 & n=2,\\
			s^3- 3 s & n=3,\\
			s^4- 4 s^2 + 2 &  n=4,\\	
			s^5-  5 s^3+ 5 s & n=5,\\
			s^6- 6 s^4 + 9 s^2-2& n=6.
		\end{cases}
	\end{equation}
	We denote by $\KK'_n$ the splitting field of rational elliptic surface
	$\Ee'_n$ over $\Q(s)$ for $1\leq n\leq 6$ which is determined  in the next sections.
	%	The singular fibers of $\Ee'_n$  are of type $II$  over  the roots of $f_n$ 	and of types $IV^*, I_0^*, IV, II, II$  at $s=\infty$, and hence applying the Shioda--Tate's formula one can see that  
	The  Mordell-weil rank  of $M'_n$ is  $2, 4, 6, 8, 8$  and  we    have $M'_n  \iso \{0\},  a_2^*, D_4^*, E_6^*, E_8, E_8, $ minimal norms  	$0, 2/3, 1, 4/3, 2, 2$,  for $n=2, \ldots, 6$ respectively.
	%	Indeed,  we    have $M'_n  \iso \{0\},  a_2^*, D_4^*, E_6^*, E_8, E_8, $	for $n=1, 2, 3, 4, 5, 6$, respectively. 
	Here, 	$a_2^*$ indicates the dual lattice of the root lattice $a_2$, etc. 
	%	The minimal norms  of these lattices are 	$0, 2/3, 1, 4/3, 2, 2$, respectively.
	%	We refer the reader to 	\cite{Shioda2008}   for   proof of  above statements as well as the 
	
	The	following  theorem is the main result of T.~Shioda on the Mordell-Weil lattice of  	$\Ee_n$.
	
	\begin{thm}
		\label{shi-thm}
		With the above notations, the invariants of  $M_n=\Ee_n (\CC(t))$ are given in Table~1, 	where $\mu_n$ denotes the length of  minimal sections.
		Moreover, the lattice $M_n $ is generated by the points   $P=(x(t), y(t) )$
		with the coordinates
		\begin{align*}
			x(t)& =\frac{a_0 + a_1 t + a_2 t^2+ a_3 t^3 + a_4 t^4}{t^2},  \ (a_i \in \CC)\\
			y(t)& =\frac{b_0 + b_1 t + b_2 t^2+ b_3 t^3 + b_4 t^4+b_5 t^5 + b_6 t^6}{t^3}, \ (b_j \in \CC) .
		\end{align*}
		More precisely, for $n=2$, a set of independent generators of $M_2$ is given by
		$(\alpha,   t  + 1/t)$ and $ (\alpha', t-1/t)$
		where  $\alpha$ and $\alpha'$ run over the roots of  cubic polynomials $u^3-2$ and $u^3+2$, respectively.
		For  $n> 2$, the lattice  $M_n$ is generated by    following   set of points:
		\begin{itemize}
			\item [(i)]	 In  cases   $n=3, 5$:
			$$   \left(x'\left(t+\frac{1}{t}\right), y'\left(t+\frac{1}{t}\right)\right) \ \text{and}\
			\left(x'\left(\zeta_n t+ \frac{1}{\zeta_n  t}\right), y'\left(\zeta_n t+ \frac{1}{\zeta_n t}\right)  \right)$$
			
			\item [(ii)]  In cases   $n=4, 6$:
			$$  \left(x'\left(t+\frac{1}{t}\right), y'\left(t+\frac{1}{t}\right)\right) \ \text{and}\
			\ \left(-x'\left(\zeta_{2n} t+\frac{1}{\zeta_{2n} t}\right) , \I \, y'\left(\zeta_{2n} t+\frac{1}{\zeta_{2n} t}\right)  \right). $$
			%$ g'_n\left( x(t+1/t), y(t+1/t), t\right) )$		 
		\end{itemize}
		where 	$(x'(s), y'(s))$  belongs to a  generating set of $M'_n$
		with the coordinates:
		$$x'(s)=a_0+ a_1 s+ a_2 s^2, \ \text{and}\  y'(s) =b_0+ b_1 s+ b_2 s^2 + b_3 s^3,  \   (a_i, b_j \in \CC).$$
		
		\begin{table}[htbp]
			\caption{Invariants of the lattices $M_n=\Ee_n(\CC(t))$}\label{Tab1}
			\begin{tabular}{|ccccccc|}
				\hline
				$n$ &  $1$ & $2$ & $3$ & $4$ & $5$ & $6$ \\ 
				\hline
				%			$M_n$ &  $\{0\}$ & $a_2^*$ & $D_4^*$ & $E_6^*$ & $E_8$ & $E_8$   \\ [5pt]
				$r_n$ &  $0$ & $4$ & $8$ & $12$ & $16$ & $16$   \\ [5pt]
				$\det(M_n)$ &  $1$ & $2^4/3^3$ & $3^4/4^2$ & $4^4/3^2$ & $5^4$ & $6^4$   \\ [5pt]
				$\mu_n$ &  - & $4/3$ & $2$ & $8/3$ & $4$ & $4$ \\ 
				\hline	 
			\end{tabular}
		\end{table}
	\end{thm}
	
	In \cite[Thm. 2.5]{Shioda2008}, Shioda proved the above theorem, but he did not determined exactly neither the coefficients nor  splitting fields $\KK_n$, which   is our main task in this paper. We refer the reader to  see the proofs of  Theorems~2.4 and 2.6 in \cite{Shioda2008} to see more  details. Here, we just provide a sketch of the main idea of the proofs.
	
	Letting $T=t^n$, $w=T+ 1/T$, and $L_n=M'_n  [2]$ for $1\leq n \leq 6$, 
	considering the elliptic $\Ee_0: y^2=x^3+w$ over $\CC(w)$, we have
	$\Ee_0 \cong \Ee_1$ and so $ M_n=\Ee_0 (\CC(t)), \ L_n=\Ee_0(\CC(s)), $ and $  N_n= \Ee_0(\CC(T)).$
	We note that $\CC(t)$  is a Galois extension of $\CC(w)$  with  Galois group
	$G=\left\langle \tau_0, \tau_n \right\rangle $ with $\tau_0: t \rightarrow 1/t$ and
	$\tau_n: t \rightarrow \zeta_n t$, where $\zeta_n$ is an $n$-th root of the unity.
	In the terminology of  Galois Theory, the fields $\CC(s)$  and $\CC(T)$ correspond
	to the subgroups $\left\langle  \tau_0  \right\rangle$ and
	$\left\langle  \tau_n \right\rangle$,  and the invariant sublattices of $M_n$  are $ L_n$ and $N_n$, respectively.
	
	By \cite[Lemma 7.2 and  7.3]{Shioda2008}, we have
	$L_n\cap N_n=\{ 0\}$, and $L_n\oplus N_n$ is an orthogonal direct sum of lattices.
	Moreover, if we let $\tilde{L}_n =\tau_n(L_n) \subseteq M_n$, then
	$\tilde{L}_n=\Ee_0 (\CC(s'))$ with $s'=\tau_n(s)= \zeta_n t+ \frac{1 }{\zeta_n t} $ such that
	$L_n \cap \tilde{L}_n =\{0\}$ for odd $n$ and
	$L_n \cap \tilde{L}_n\cong M'_2 $  otherwise.
	In \cite[Lemma~7.4]{Shioda2008}, it is proved that $M_n=L_n + \tilde{L}_n$ for $n=3,5$ and $\det(M_n)$
	is equal to $3^4/4^2$ for $n=3$, and $5^4$ for $n =5$.
	In the case of $n=4, 6$,    denoting the fourth root of the unity  by $\I$,
	redefining $\tilde{L}_n$ as the image of $L_n$ by   the  following automorphism of  $M_n$,
	\begin{equation}
		\label{phimap}
		\phi_n: (x(t), y(t))\rightarrow (-x( \zeta_{2n}  t),\ \I \, y (\zeta_{2n}  t)),
	\end{equation}
	and using \cite[Lemma~7.5]{Shioda2008}, we have  $L_n \cap \tilde{L}_n =\{0\}$
	and $\det(L_n + \tilde{L}_n)= 4^4/3^2$ for $n=4$ and $6^4$ for $n=6$.
	Therefore, one may conclude that
	$ N_n\oplus L_n \oplus \tilde{L}_n$  is a sublattice of finite index in $M_n$ for $n=4, 6$.  
	
	%This is a restatement of \cite[Thm. 2.4, 2.6]{Shioda2008}

	%++++++++++++++++++++++++++++++++++++++++++++
	\section{An algorithmic approach to   the proof of the theorems}
	%++++++++++++++++++++++++++++++++++++++++++++	
	In this section,  we provide an algorithmic approach for   proof of all results of the paper.
	By  Shioda's results \ref{shi-thm}, to determine the splitting field $\KK_n$ of $\Ee_n$ and a set of the linearly independent generating points of $\Ee_n ( \KK_n )$, we  will do the  steps provided in Table~\ref{Tab3}
	
	\begin{center}
		\label{Tab3}
		\begin{table}[h]
			\caption{Algorithm for computation on $\Ee_n(\KK_n)$}
			\begin{tabular}{|c l |}
				\hline  
				& Computing the Splitting field and Generators of  $\Ee_n (\KK_n)$   \\
				\hline
				{\bf Input:}  &  Defining equation  of elliptic curve $\Ee_n$  over $\Q(t)$ of rank $r_n$\\
				& over $\CC(t)$ with  known invariant as in Table \ref{Tab1}  \\
				\hline 	
				{\bf Step 1:} &  Determining the splitting field and linearly generators  of rational\\
				&  elliptic surface   $\Ee'_n:  y^2=x^3+ f_n(s)$ over $\Q(s)$ of rank $r'_n$ over $\CC(s)$\\
				&   $\bullet$   Take   points (sections) of  elliptic surface  of the form \\
				& $(x'(s), y'(s))=(a_0+ a_1 s+a_2 s^2, b_0+ b_1 s+ b_2 s^2 + b_3 s^3),$   and \\
				& substitute  into  the  	equation of $\Ee'_n$, to get a set of equations \\
				& in $a_i's$ and $b_j's$  defining an ideal in $\Q[a_0, a_1, a_2, b_0, b_1, b_2, b_3]$\\
				&  $\bullet$ Finding the fundamental  polynomial   of above ideals using \\
				&  the  command {\sf UnivariatePolynomial}  of package {\sf PolynomialIdeals}\\
				&  in  {\sf Maple} and factoring it to linear factors, as given in \cite{k3-codes}  \\
				&   $\bullet$ Use {\sf Pari/GP} code  in \cite{k3-codes} to find a defining minimal\\ 
				& polynomial $g'_n(x)$ of  the splitting field of the fundamental polynomials, \\
				&  i.e, defining minimal polynomial of the splitting  fields $\KK'_n$\\
				&  $\bullet$ Choose a set of appropriate roots of fundamental polynomials\\
				&   to get   linearly independent generators of  $\Ee'_n(\KK'_n)$\\
				\hline 					  
				{\bf Step 2:}    &      Determining the splitting field $\KK_n$ and linearly generators of $\Ee_n(\KK_n)$ \\
				&  $\bullet$ Use {\sf  Pari/GP} and {\sf SageMath} k3-codes to find a defining minimal\\ 
				& polynomial $g_n(x)$ of compositum field $\KK_n=\KK'_n(\zeta_m)$  with  \\
				& $m=n$  for  $n=3,5$ and $m=2n$ for $n=4,$ and $6$\\	                  
				& $\bullet$ Transforming the points $(x'(s), y'(s))\in \Ee'_n(\KK'_n)$
				into points\\
				&  belonging   $\Ee_n(\KK_n)$ using the transformations given in \ref{shi-thm}\\
				\hline 
				{\bf Output:} & The splitting field $\KK_n \subset \CC$ of  $\Ee_n$ and  a set of  \\
				& linearly independent generators  for  $\Ee_n(\KK_n)$ \\	
				\hline 
			\end{tabular}
		\end{table}
	\end{center}

	%++++++++++++++++++++++++++++++++++++++++++++
	\section{The cases    $ \Ee_2 $ and $\Ee_3$}
	%++++++++++++++++++++++++++++++++++++++++++++
	\label{case123}
	In this section, we consider the Mordell--Weil lattices of the  simple cases $\Ee_2, $ and $\Ee_3$.
	\subsection{Proof of Theorem~\ref{main-0}}
	The structure of Mordell--Weil lattice of    $\Ee_2$ over $\CC(t)$  is treated in \cite[Theorems 6.1]{Shioda2007b} and 
	\cite[Theorem.~7.1]{Shioda2008}. 
	
	In the case of $\Ee_2$, by Theorem~\ref{shi-thm}, a set of independent generators
	can be found between   points of the form 	$\left(a, b t + c+ d/t \right).$ 
	Substituting these points in the equation of    $\Ee_2$ leads to  $c=0$,  $b, d \in \{ \pm 1\}$. If $b$  and $d$ have the same sign, then $a^3-2=0$ and otherwise $a^3+2=0$. Hence, there are totally six points and 
	one can check that the Gram matrix of the points  $P_1, P_2, P_3, P_4$  given in the statement of  Theorem~\ref{main-0} is
	\begin{equation}
		\label{R2}
		R_2=\frac{2}{3}
		\begin{pmatrix}
			2 & 1 & 0 & 0
			\\ \noalign{\medskip}
			1 & 2 & 0 & 0
			\\ \noalign{\medskip}
			0 & 0 & 2 & 1
			\\ \noalign{\medskip}
			0 & 0 & 1 & 2
		\end{pmatrix},
	\end{equation}
	which  has the determinant $2^4/3^2$ as desired. 
	The splitting field  $\KK_2$ is equal to an extension of $\Q$ with contains the roots of 
	$a^3-2=0$ and otherwise $a^3+2=0$, say $\KK_2=\Q(\zeta_3, 2^{1/3})$ with a minimal  defining polynomial $x^6+ 108.$

	%Thus the proof of  Theorem~\ref{main-0} is completed. 
	%--------------------------------
	\subsection{Proof of Theorem~\ref{main-a}}
	\label{coef3}
	
	We consider the rational elliptic surface $\Ee'_3: y^2=x^3-(s^3-3s)$ with 	 discriminant  $27 s^2 (s^2-3)^2.$
	%	The singular fibers of   $\Ee'_3$ are  of type $II$	over $s=0, \pm \sqrt{3}$, and   of type $I_0^*$   over $s=\infty$.	Applying the Shioda--Tate's formula shows that
	According to Shioda's result,	the rank of  $\Ee'_3(\CC(s))$ is equal to $4$ and $\Ee'_3(\CC(s)) \iso D_4^*$.
	To find a set of independent generators, we consider  the points  $Q=(a s+ b, c s + d)$ and 
	substitute its   coordinates in the equation of $\Ee'_3: y^2=x^3-(s^3-3s)$ to obtain the following  equalities:
	\begin{equation}
		\label{IDEQ3}	
		\begin{aligned}[b] a^3+1&=0, & c^2-3 a^2 b &=0, & -3 ab^2 + 2 cd +3&=0,  & d^2-b^3&=0.\end{aligned}
	\end{equation}	
	Form the second and third equities, we get 
	\begin{equation}
		\label{IDEQ3-1}	
		\begin{aligned}[b] b& =c^2/3 a^2, & d&=-(c^4+1)/6c.\end{aligned}
	\end{equation}
	Hence, the last equality  gives us $c^8-54 c^4- 243=0, $  
	whose 	roots are as follows:
	
	%	$$c=\pm \sqrt{3}	(3+ 2 \sqrt{3})^{\frac{1}{4}}, \    \I \pm \sqrt{3}	(3+ 2 \sqrt{3})^{\frac{1}{4}}, $$
	%	$$\ \ \ \pm  \frac{\sqrt{2}(1+\I)}{2} (3- 2 \sqrt{3})^{\frac{1}{4}},	    \pm  \frac{\sqrt{2}(1-\I)}{2} (3- 2 \sqrt{3})^{\frac{1}{4}},$$	
	
	%	$$c=\pm 3^{\frac{5}{8}}	\epsilon_1^{\frac{1}{4}}, \  \pm  \I 3^{\frac{5}{8}}	\epsilon_1^{\frac{1}{4}}, \  \pm 3^{\frac{5}{8}}  \zeta_8 {\epsilon'_1}^{\frac{1}{4}}, \    \pm 3^{\frac{5}{8}}\zeta'_8 {\epsilon'_1}^{\frac{1}{4}},$$
	
	\begin{equation}
		\label{eps1-2}	
		\begin{aligned}[b] 
			c=&\pm \sqrt{3}	(3+ 2 \sqrt{3})^{\frac{1}{4}}, \ \ \pm  \frac{\sqrt{2}(1+\I)}{2} (3- 2 \sqrt{3})^{\frac{1}{4}},\\
			& \pm \I  \sqrt{3}	(3+ 2 \sqrt{3})^{\frac{1}{4}}, \ \ 
			\pm  \frac{\sqrt{2}(1-\I)}{2} (3- 2 \sqrt{3})^{\frac{1}{4}}.
		\end{aligned}
	\end{equation}

	%	$$\pm\, \sqrt{3}\, \alpha_1,  \
	% \ -\beta,     \I \, \beta,
	%	\pm \, \I \, \sqrt{3} \, \alpha_1, \ \pm \alpha_2   \frac{\sqrt{6} \, (1+\I)}{2}, \  	%\alpha \cdot \left( \frac{-1+\I}{2}\right), \
	%	\pm \alpha_2 \frac{\sqrt{6} \, (1-\I)}{2} , \	%\alpha \cdot \left( \frac{1-\I}{2}\right),
	%	\text{with} \
	%	\alpha_1 = (2 \sqrt{3} +3)^{\frac 1 4}, \ \text{and} \
	%	\alpha_2 = (2 \sqrt{3} -3)^{\frac 1 4}.$$
	%where
	%$$\displaystyle \alpha= (2 \sqrt{3} -3)^{\frac 1 4}, \ \text{and} \ \beta = (2 \sqrt{3} +3)^{\frac 1 4},$$
	The above eights roots  together with  the three roots of $a^3+1=0$, say $a=-1, (1\pm \I\sqrt{3})/2$,  determine $24$ points
	on $\Ee'_3$ generating  Mordell--Weil lattice
	$\Ee'_3 (\CC(s))$.
	The points with $a=-1$ generate a sublattice isomorphic to the unit matrix of degree four. 
	By  straight  computations and similar argument as  in
	\cite[Section~6]{Shioda1991d}, one can check that 
	four points $Q_j=(a_j s+b_j, c_j s+d_j)$ generate  $\Ee'_3 (\CC(s))$,  
	where  their coefficients are $$a_1=a_2=a_3=-1, a_4=\frac{1+\I \sqrt{3}}{2}$$ and
	$$	\begin{aligned}[b] 
		b_1&=\sqrt{3+2 \sqrt{3}}, & c_1&=\sqrt{3}\, \left(3+2 \sqrt{3}\right)^{\frac{1}{4}}, &d_1&
		=-\left(3+2 \sqrt{3}\right)^{\frac{3}{4}},\\
		b_2&= -\sqrt{3+2 \sqrt{3}},& c_2& =\mathrm{I} \sqrt{3}\, \left(3+2 \sqrt{3}\right)^{\frac{1}{4}}, 
		& d_2&=\mathrm{I} \left(3+2 \sqrt{3}\right)^{\frac{3}{4}},\\
		b_3&=\mathrm{I} \sqrt{2 \sqrt{3} -3}, & c_3& =\frac{ \sqrt{6}(\mathrm{I}+1)}{2} \left(2\sqrt{3}-3 \right)^{\frac{1}{4}}, 
		& d_3&=\frac{ \sqrt{2}(\mathrm{I}-1)}{2} \left(2 \sqrt{3}-3\right)^{\frac{3}{4}},\\
		b_4&=-\frac{ \left(1+\mathrm{I} \sqrt{3}\right)}{2}\sqrt{3+2 \sqrt{3}}, & c_4&=\sqrt{3}\, \left(3+2 \sqrt{3}\right)^{\frac{1}{4}}, 
		&d_4&=-\left(3+2 \sqrt{3}\right)^{\frac{3}{4}}.
	\end{aligned}$$

		The Gram matrix of the  points $Q_j$'s has  determinant $1/4$ and is given by
		\[R'_3=\frac{1}{2}
		\begin{pmatrix}
			2 & 0 & 0 & 1
			\\ \noalign{\medskip}
			0 & 2 & 0 & 1
			\\ \noalign{\medskip}
			0 & 0 & 2 & 1
			\\ \noalign{\medskip}
			1 & 1 & 1 & 2
		\end{pmatrix}.
		\]
		Thus, the splitting field $\KK'_3$  of $\Ee'_3$ over $\Q(t)$ is compositum of the fields defined by the polynomials  $a^3+1=0$ and
		$c^8-54 c^4- 243=0$, containing the field $\Q(\zeta_3, (3+ 3 \sqrt{3})^{\frac{1}{4}})$.
		% Hence  $\KK'_3$ must contains 	$\Q \left( \zeta_6,   3^{\frac{1}{8}}, 
		%	\epsilon_1^{\frac{1}{4}},  \epsilon_2^{\frac{1}{4}} \right)$ including
		% all of their roots.

		Using Theorem~\ref{shi-thm},
		and  substituting  $s=t + 1/t$  and
		$s=\zeta_3 t + \frac{1}{\zeta_3 t}$
		in the coordinates of  $Q_j$'s for $j=1, 2,3, 4$, we  obtain
		$$\begin{aligned}
			P_j  &= (x(t), y(t)) =\left( \frac{a_j t^{2}+b_j t +a_j}{t}, \ \frac{c_j t^{2}+d_j t +c_j}{t}\right),\\
			P_{j+4} &= (x(\zeta_3 t), y(\zeta_3 t)) = \left(\frac{a_j \zeta_3^{2} t^{2} +b_j \zeta_3 t +a_j}{\zeta_3 t}, \
			\frac{c_j \zeta_3^{2} t^{2}+d_j \zeta_3 t +c_j}{\zeta_3 t} \right).
		\end{aligned}$$
		By the properties of  height pairing and knowing that  $\KK_3(t)$ is a quadratic extension of $\KK_3(s)$,
		where $\KK_3=\KK'_3(\zeta_3)=\KK'_3$.

		Using Sagemath, we   obtained the minimal defining polynomial of $\KK_3$,  containing $\Q(\zeta_6, (3+ 3 \sqrt{3})^{\frac{1}{4}})$,  is  equal to  the compositum of   fields defined by  $x^3-1=0$,  $x^3+1=0$ and $x^8-54 x^4- 243=0$ having  a minimal defining polynomial
		as follows:
		\begin{align}
			g_3(x) &= x^{16} + 8x^{15} + 36x^{14} + 112x^{13} + 158x^{12} - 144x^{11} - 836x^{10}  \notag \\
			&\quad  + 86040x^{5} - 1144x^{9} + 3051x^{8}  + 14624x^{7} + 45820x^{6}\notag \\
			&\quad  + 109130x^{4}  + 91912x^{3} - 60552x^{2} - 94600x + 49141.
			\label{f3}
		\end{align}
		%which can be factored completely into linear factors  over $\Q(\I, (3+ 3 \sqrt{3})^{\frac{1}{4}})$, see \cite{k3-codes}.

		By properties of height pairing,  we have
		$$\left\langle P_i, P_{i+j} \right\rangle  = - \frac{1}{2}  \left\langle P_i, P_j \right\rangle \  (1\leq i, j \leq 4).$$
		Using this fact and the  matrix $R'_3$,
		one can see that the Gram matrix of the eight points $P_1, \ldots, P_8$ is
		\[R_3=\frac{1}{4}
		\begin{pmatrix}
			4 & 0 & 0 & 2  & -2 & 0 & 0 & -1 \\  \noalign{\medskip}
			0 & 4 & 0 & 2  & 0 & -2 & 0 &  -1 \\   \noalign{\medskip}
			0 & 0 & 4 & 2  & 0 & 0 & -2 &  -1 \\  \noalign{\medskip}
			2 & 2 & 2 & 4  & -1 & -1 & -1 & -2\\  \noalign{\medskip}
			-2 & 0 & 0 & -1  &  4 & 0 & 0 & 2\\  \noalign{\medskip}
			0 & -2 & 0 &  -1  &  0 & 4 & 0 & 2 \\  \noalign{\medskip}
			0 & 0 & -2 &  -1   &  0 & 0 & 4 & 2 \\  \noalign{\medskip}
			-1 & -1 & -1 & -2 &	2 & 2 & 2 & 4
		\end{pmatrix}.\]
		and its determinant is $3^4/4^2$ as given by Theorem~\ref{shi-thm}.
		%	Therefore, we have finished the proof of  Theorem~\ref{main-a}.
		We refer the reader to see \cite[check-3]{k3-codes} for the computations of this section.

		%++++++++++++++++++++++++++++++++++++++++++++
		\section{The case of  $\Ee_4$}
		%++++++++++++++++++++++++++++++++++++++++++++4
		\label{case4}
		In this section,  we prove Theorem~\ref{main-b} using   the following result on the rational elliptic surface
		$\Ee'_4.$

		\begin{thm}
			\label{main-b1}
			The splitting field  $\KK'_4$ of rational elliptic surface
			$$ \Ee'_4 : y^2=x^3-(s^4- 4 s^2 + 2),$$  is   the number field $\KK'_4 $ defined by a polynomail of degree 24 given by \ref{g14}.
			%		$\Q \left( \zeta_{12}, 2^{\frac{1}{12}},  \epsilon_2^{\frac{1}{6}},  \epsilon_3^{\frac{1}{6}}\right)$
			%	with
			%	$$\begin{aligned}[b]\alpha_0 &=  2^{\frac{1}{3}}  ( 18 \sqrt{6} +44)^{\frac{1}{6}}, &
				%	\beta_0 &= 2^{\frac{1}{3}}  (4+\I \,10  \sqrt{2})^{\frac{1}{6}},\end{aligned}$$
			Moreover,   the Mordell--Weil lattice  $\Ee'_4(\KK'_4 (s))$ is generated by the   points
			\begin{align*}
				Q_j&=\left(a_j s +b_j,  s^2 + c_j s + d_j \right),
			\end{align*}
			for $j=1, \ldots, 6$, where  $a_j, b_j, c_j, d_j$ are given in Subsection~\ref{coef4}.
		\end{thm}

		\subsection{Proof of Theorem~\ref{main-b1}}
		\label{coef4}
		Since the discriminant  of $\Ee'_4$ is   $-27(s^4- 4 s^2 + 2)^2$,
		the singular fibers of $\Ee'_4$ are of type $II$
		over the roots of $s^4- 4 s^2 + 2$   and  of type  $IV$   over $s=\infty$.
		Then,  the Shioda--Tate's formula shows that  the Mordell--Weil rank of  $\Ee'_4(\CC(s)) $ is equal to $6$ and
		$\Ee'_4(\CC(s)) \iso E_6^*.$
		Based on \cite[Theorem~10.5]{Shioda1991d}, a set of six independent generators of	$\Ee'_4(\CC(s))$ 
		can be found between  the set of $27$ rational  points  $Q=(a s+ b, \  s^2 + c s + d)$.
		Substituting these in the equation of $\Ee'_4$ leads to the following equalities:
		\begin{align}
			2 c -a^{3}&=0,&
			3 a^{2} b -c^{2}-2 d -4&=0, \notag \\
			3 a b^{2}-2 c d &=0,&
			b^{3}-d^{2}+2 &=0. \label{eqs4}
		\end{align} 
		
		From the first two   qualities, we get
		\begin{equation}
			\label{4eq1}
			\begin{aligned}[b]
				c &= \frac{a^3}{2}, & d &=-\frac{1}{8} (a^{6}-12 a^{2} b +16),
			\end{aligned}
		\end{equation}
		and  two polynomials in $b$ with coefficients in the ring $\Q[a]$ as
		$$\begin{aligned}
			p_1&=b^{2} -12 a^{5} b +16 a^{3}+24 a a^{9},\\
			p_2&=-64 b^{3}+144 a^{4} b^{2}-24 a^{2} \left(a^{6}+16\right) b +a^{12}+32 a^{6}+128.
		\end{aligned}$$
		Taking the resultant  respect to $b$  of $p_1$ and $p_2$ 
		gives a polynomial of degree $27$ of the form  $\Phi(a)=a^{3}  \Phi_1(a) \Phi_2(a),$     where
		\begin{equation*}
			%\label{4eq2}
			\Phi_1(a)=\left(a^{12}-352 a^{6}-128\right), \ \text{and}   \ 
			\Phi_2(a) = \left(a^{12}-32 a^{6}+3456\right).
		\end{equation*}
		By  \eqref{4eq1} and  using  the roots of $\Phi(a)$, we obtain  coefficients of $27$
		points  $Q=(a s + b, s^2 + c s + d)$ in $ \Ee'_4 (\CC(t))$.
		The  factors of degree $12$ of $\Phi(a)$ can be decomposed as follows:
		$$\begin{aligned}[b]
			\Phi_1(a)&= \prod_{\ell=0}^5 \left(a-   
			2^{\frac{7}{12}} \zeta_{12}^{2 \ell}  \epsilon_2^{\frac{1}{6}} \right)
			\prod_{\ell=0}^5 \left( a- 2^{\frac{7}{12}} \zeta_{12}^{2 \ell +1} {\epsilon'_2}^{\frac{1}{6}} \right), \text{and} \\
			\Phi_2(a)&=  \prod_{\ell=0}^5 \left(a- 2^{\frac{7}{12}}\zeta_{12}^{2 \ell} \epsilon_3^{\frac{1}{6}}\right)
			\prod_{\ell=0}^5 \left( a- 2^{\frac{7}{12}} \zeta_{12}^{2 \ell +1} {\epsilon'_3}^{\frac{1}{6}} \right),
		\end{aligned}$$
		where $\epsilon_2, \epsilon'_2, \epsilon_3, $ and  $ \epsilon'_3$ are as in Table~\ref{tab:notation}.
		Thus, the splitting field  $\KK'_4$ of  $\Ee'_4$ is equal to compositum of the splitting field of the polynomials  $\Phi_1(a)$  and $\Phi_2(a)$, which has a minimal  defining   polynomial as:
		%\begin{equation}
		
		\begin{align}
			g_4(x)	& = x^{24} - 12x^{22} + 114x^{20} - 664x^{18} 
			+ 2856x^{16} - 8928x^{14} + 21196x^{12} \notag \\
			& \quad - 33576x^{10} + 35484x^{8} - 20544x^{6} + 5832x^{4} - 720x^{2} + 36. 
			\label{g14}	 
		\end{align} 
		%	\end{equation}
	
	%Note that 	$\KK'_4 $ contains  $ \Q \left( \zeta_{12},  2^{\frac{1}{12}}, \epsilon_2^{\frac{1}{6}}, \epsilon_3^{\frac{1}{6}} \right)$,	for which the  polynomial $\Phi(a)$	can be written as a product of  linear factors.
	Thus,  $\Ee'_4 (\CC(s)) =\Ee'_4 (\KK'_4(s))$ and
	by  straight height computations and  the determinant of lattice $\Ee'_4 (\CC(s)) \iso E_6^*$, we obtained its six independent generators  $Q_j=(a_j s + b_j, s^2 + c_j s + d_j)$ 
	with the coefficients as follows: 
	$$	\begin{aligned}[b]
		a_1 & = 0,  &  b_1 = & 2^{\frac{1}{3}}, & c_1 & =0, & d_1 &= -2,\\
		a_2 & =  2^{\frac{7}{12}} \epsilon_2^{\frac{1}{6}},  &  
		b_2 = & 2^{\frac{5}{6}}  (\sqrt{2} + \sqrt{3}), & 
		c_2 & =2^{\frac{3}{4}} \epsilon_2^{\frac{1}{2}}, & 
		d_2 &=   3 \sqrt{2}  (\sqrt{2} + \sqrt{3}),\\
		a_3 & = 2^{\frac{7}{12}} \epsilon_3^{\frac{1}{6}},  &  
		b_3 = &  \frac{2^{\frac{1}{3}} \epsilon_3^{\frac{2}{3}} (\I \epsilon_3 +1 )}{9},
		%	\left( \frac{4- \I \sqrt{2}}{9}\right), 
		& c_3 & = 2^{\frac{3}{4}} \epsilon_3^{\frac{1}{2}},
		& d_3 &=  2  +\I \sqrt{2},\\
		a_4 & = 0,  &  b_4 = & \zeta_{3}   2^{\frac{1}{3}}, & c_4 & =0, & d_4 &=-2,\\
		a_5 & = 2^{\frac{7}{12}} \zeta_{12}   {\epsilon'_2}^{\frac{1}{6}},& 
		b_5 =& 2^{\frac{5}{6}} \zeta_{6}^5  {\epsilon'_2}^{\frac{2}{3}}(\sqrt{2} + \sqrt{3}),&
		c_5 & = \I 2^{\frac{3}{4}} {\epsilon'_2}^{\frac{1}{2}}, & 
		d_5 &= 3 \sqrt{2}  (\sqrt{2} - \sqrt{3}),\\
		a_6 & =   \zeta_{6} 2^{\frac{7}{12}} \epsilon_3^{\frac{1}{6}},  & 
		b_6 = &  \frac{- 2^{\frac{1}{3}} \epsilon_3^{\frac{2}{3}} (\epsilon_3-\I) }{9 \zeta_{12}}, 
		%\left( 4 \zeta_{6} + \frac{ \sqrt{2}}{ \zeta_{12}} \right), 
		& c_6 & =- 2^{\frac{3}{4}} \epsilon_3^{\frac{1}{2}},
		& d_6 &=  2  +\I \sqrt{2}.
	\end{aligned}$$
	The Gram matrix of the points $Q_1, \ldots, Q_6$ is given by
	\[R'_4= \frac{1}{3}
	\begin{pmatrix}
		4 & -2 & 1 & -2 & 1 & -2 \\ \noalign{\medskip}
		-2 & 4 & 1 & 1 & 1 & 1 \\ \noalign{\medskip}
		1 & 1 & 4 & -2 & 1 & 1 \\ \noalign{\medskip}
		-2 & 1 & -2 & 4 & -2 & 1 \\ \noalign{\medskip}
		1 & 1 & 1 & -2 & 4 & -2 \\ \noalign{\medskip}
		-2 & 1 & 1 & 1 & -2 & 4
	\end{pmatrix},\]
	which is of determinant $1/3$ as desired. Therefore, they are independent generators of $\Ee'_4(\KK'_4 (s))$.

	\subsection{Proof of Theorem~\ref{main-b}}
	\label{m-b}
	Considering  Theorem~\ref{shi-thm}
	and  substituting  $s=t + 1/t$
	in the coordinates of  points  $ Q_j= \left(  a_j s+ b_j,  s^2 + c_j s + d_j \right) \in  \Ee'_4(\KK'_4(s))$,
	we  obtain
	$$P_j=\left( \frac{a_j t^{2}+b_j t +a_j}{t},  \frac{t^{4}+ c_j t^{3}+(d_j +2) t^{2}+c_j t +1}{t^{2}} \right),$$
	and    their images   $P_{j+6} = \phi_4(P_j)$, under the automorphism $\phi_4$ of  $\Ee_4$, with coordinates
	%	$$\phi_n: (x(t), y(t))\rightarrow (-x( \zeta_{2n}  t),\ \I \, y (\zeta_{2n}  t)),$$ 
	%	$$g'_4: \left( x\left(t+ \frac{1}{t}\right), y\left(t+ \frac{1}{t}\right)\right) \rightarrow
	%	\left(-x\left(\zeta_{8}\, t+ \frac{1}{\zeta_{8}\, t}\right), \I \, y\left(\zeta_{8}\, t+ \frac{1}{\zeta_{8}\, t}\right)  \right),$$
	%	we get  the following points
	\begin{align*}
		x(P_{j+6}) & =- \frac{\left(1+\mathrm{i}\right) a_j t^{2}+2 b_j t +\left(2-2 \mathrm{i}\right) a_j}{2 t}, \\ 
		y(P_{j+6}) &= \frac{ \mathrm{i}  t^{4}+\left(1+\mathrm{i}\right) c_j t^{3}+\left(4+2 d_j \right) t^{2}+\left(2-2 \mathrm{i}\right) c_j t -4 \mathrm{i}}{2 t^{2}}, 
	\end{align*}
	for $j=1, \ldots, 6$, which all  together generates $\Ee_4 (\CC(t))=\Ee_4 (\KK_4(t))$, where $\KK_4=\KK'_4(\zeta_8)=\KK'_4$, because  the compositum of the polynomials $g_4(x)$ and $x^8-1$ is  leads to the same number field.
	The Gram matrix of the points $P_1, \ldots, P_{12}\in \Ee_4(\KK_4(t))$ is given by
	{	\small 	
		$$R_4= \frac{1}{3}
		\left(\begin{array}{cccccccccccc}
			8 &-4 & 2 &-4 & 2 &-4 & 0 & 0 & 0 & 0 & 0 & 0 \\ \noalign{\medskip}
			-4 & 8 & 2 & 2 & 2 & 2 & 0 & 0 & 0 & 0 & 0 & 0 \\ \noalign{\medskip}
			2 & 2 & 8 &-4 & 2 & 2 & 0 & 0 & 0 &-3 & 0 & 0 \\ \noalign{\medskip}
			-4 & 2 &-4 & 8 &-4 & 2 & 0 & 0 & 3 & 0 & 0 & 0 \\ \noalign{\medskip}
			2 & 2 & 2 &-4 & 8 &-4 & 0 & 0 & 0 & 0 & 0 & 0 \\ \noalign{\medskip}
			-4 & 2 & 2 & 2 &-4 & 8 & 0 & 0 & 0 & 0 & 0 & 0 \\ \noalign{\medskip}
			0 & 0 & 0 & 0 & 0 & 0 & 8 &-4 & 2 &-4 & 2 &-4 \\ \noalign{\medskip}
			0 & 0 & 0 & 0 & 0 & 0 &-4 & 8 & 2 & 2 & 2 & 2 \\ \noalign{\medskip}
			0 & 0 & 0 & 3 & 0 & 0 & 2 & 2 & 8 &-4 & 2 & 2 \\  \noalign{\medskip}
			0 & 0 &-3 & 0 & 0 & 0 &-4 & 2 &-4 & 8 &-4 & 2 \\ \noalign{\medskip}
			0 & 0 & 0 & 0 & 0 & 0 & 2 & 2 & 2 &-4 & 8 &-4 \\ \noalign{\medskip}
			0 & 0 & 0 & 0 & 0 & 0 &-4 & 1 & 2 & 2 &-4 & 8
		\end{array}\right),$$}
	and its determinant is $4^4/3^2$ as desired. Therefore, the proof of Theorem~\ref{main-b} is completed.
	We refer the reader to see \cite[check-4]{k3-codes} for the computations of this section.

	%++++++++++++++++++++++++++++++++++++++++++++
	\section{The case of  $\Ee_5$}
	%++++++++++++++++++++++++++++++++++++++++++++
	\label{case5}
	In this section, we  prove  Theorem~\ref{main-c} using the 
	following result on the splitting  field and   a set of independent generators  of
	$$\Ee'_5:  y^2=x^3 +  s^5 - 5 s^3 + 5 s,$$ 
	over $\CC(s)$.  
	%	The reader can see \cite[check-5]{k3-codes} for the computations of this section.
	
	%The reader can compare this result with \cite[Corollary~6]{Shioda2006a}.
	
	\begin{thm}
		\label{main-c1}
		The splitting field  $\KK'_5$ of  $ \Ee'_5$  is defined by a polynomial of degree 96 given  \cite{k3-codes},  
		which contains
		$\Q \left( \zeta_{12},  5^{\frac{1}{24}},  (\epsilon_4 \epsilon_5)^{\frac{1}{2}}\right)$.
		%			$$	\zeta_5=\frac{\sqrt{5}-1 + \I \sqrt{2} \sqrt{5 +\sqrt{5}}}{4},\$$$	
		Moreover, the lattice   $\Ee'_5(\KK'_5 (s))$ is generated by the   points
		\begin{align*}
			Q_j=(x_j(s), y_j(s))&=\left(\frac{s^2 +a_j s +b_j}{u_j^2}, \frac{s^3 + c_j s^2 + d_j s + e_j}{u_j^3} \right),
		\end{align*}
		for $j=1,\ldots, 8$, where  $a_j, b_j, c_j, d_j, e_j$ are given  in \cite[check-5]{k3-codes}, and  
		the constants $u_j$'s are as follows:
		$$\begin{aligned}[b]
			u_1 & = \I   5^{\frac{1}{24}}  (\epsilon_4 \epsilon_5^{-1})^{\frac 12},&
			u_2 & = \I   5^{\frac{1}{24}}  (\epsilon_4^{-1} \epsilon_5)^{\frac 12},&
			u_3 & = \I   5^{\frac{1}{24}}  (\epsilon_4  \epsilon_5)^{\frac 12},\\
			u_4 & = \I   5^{\frac{1}{24}}   (\epsilon_4  \epsilon_5)^{-\frac 12},&
			u_5 & = \I   5^{\frac{1}{24}}   (\epsilon_4\epsilon_5^{-1})^{\frac 12}  \epsilon_6,  &
			u_6  &= \I   5^{\frac{1}{24}}   (\epsilon_4^{-1}  \epsilon_5)^{\frac 12}  \epsilon_6, \\
			u_7 & = \I   5^{\frac{1}{24}}   (\epsilon_4 \epsilon_5 )^{\frac 12} \epsilon_6^{-1}, &
			u_8 & = \I   5^{\frac{1}{24}}   (\epsilon_4\epsilon_5  \epsilon_6 )^{-\frac 12}, &
		\end{aligned}$$
		where  $\epsilon_4,\epsilon_5$, and $\epsilon_6$ are as in th estatement of Theorem~\ref{main-c}.
		%		$$	\zeta_{12}= e^{2\pi i/12}=\frac{\I + \sqrt {3}}{2}, $$
		%			$$\begin{aligned}[b]\epsilon_4 &=1-\zeta_{12}, & \epsilon_5&=\frac{1 + \sqrt {5}}{2} \zeta_{12}, &
			%				\epsilon_6 &=\frac{\sqrt {3}-\sqrt {5}}{2}(\zeta_{12}+ \zeta_{12}^{10}),\end{aligned}$$
		%			are the fundamental units of the number field $\Q(\I,\sqrt {3},\sqrt {5})=\Q(\zeta_{12},\sqrt{5})$.
	\end{thm}
	
	%---------------------------------------------------------
	\subsection{Proof of Theorem~\ref{main-c1}}
	\label{fpe0}
	
	%We consider the elliptic surface  $\Ee'_5:  y^2=x^3 + (s^5 - 5 s^3 + 5 s)$.
	Since  $\Ee'_5(\CC(s))$ is isomorphic to $E_8$,  
	%	by the general theory \cite{Shioda1990a} of rational elliptic surfaces of rank $8$, 
	there are  $240$ points $Q \in \Ee'_5(\CC(s))$, corresponding to the $240$ minimal roots of $E_8$,   of the form:
	\begin{equation*}
		%\label{eq3}
		Q=\left(\frac{s^2 +a s +b}{u^2}, \frac{s^3 + c s^2 + d s + e}{u^3} \right),
	\end{equation*}
	for suitable constants $a, b, c, d, e, u \in \CC.$
	Substituting the coordinates of $Q$ in the equation   of $\Ee'_5$, we get  the following six relations:

	\begin{align}
		2c -3a -u^6  &=0,  & 
		2d-3a^2+c^2-3b  &=0, \notag \\
		2e-a^3-6ab+2cd+5 u^6 & =0, & 
		3b^2 +3a^2b - 2ce -d^2 & =0, \notag \\
		3a b^2+5 u^6 -2de  &=0, &	
		b^3-e^2  &= 0. 	 	\label{eq4a} 	
	\end{align} 

	By the first three relations, we obtain $c, d, e$ in terms of $a, b, $ and $u$ as:
	\begin{equation}
		\label{eq4}
		\begin{aligned}[b]
			c &= \frac{3a +u^6}{2}, &
			d &=\frac{3a^2-c^2+3b}{2}, &
			e  &=\frac{a^3+6ab-2cd-5 u^6}{2}.
		\end{aligned}
	\end{equation}
	Using Maple,  we calculate the fundamental polynomial of the Ideal generated by equations~\ref{eq4a} of degree 240 in variable $u$, see \cite[check-5]{k3-codes}.

	Letting $V=u^{12}$, the polynomial $\Phi (V)$
	decomposes into  four irreducible factors in $\Z[V]$ over $\Q$, namely,	
	$$\begin{aligned}
		\Phi_1(V) &=  V^4-56700 V^3-1204210 V^2-283500 V+25, \\
		\Phi_2(V) &=  V^4+6660 V^3-685810 V^2-91320300 V+25,\\
		\Phi_3(V) &=  V^4-1260 V^3+1178590 V^2-4592700 V+13286025,\\
		\Phi_4(V)& =V^8-24300 V^7+280019230 V^6\\
		&\quad -18498253500 V^5 +569262158025 V^4 +5919441120000 V^3  \\
		&\quad +28673969152000 V^2 +796262400000 V+10485760000.
	\end{aligned}$$
	%\end{equation}
	%and one more $F_4(V)$ of degree $8$.
	Then,  one can decompose all polynomials $\Phi_1$, $\Phi_2$,  $\Phi_3$ and $\Phi_4$ over
	$k_0=\Q(\I, \sqrt{3}, \sqrt{5})=\Q(\zeta_{12}, \sqrt{5})$ into the product of linear factors as given below,
	$$	\begin{aligned}
		\Phi_1(V) &=  (V-v_1)(V-v_1^\sigma)(V-v_1^\tau)(V-v_1^{\sigma \tau}), \\
		\Phi_2(V) &=  (V-v_2)(V-v_2^\sigma)(V-v_2^\tau)(V-v_2^{\sigma \tau}),\\
		\Phi_3(V) &=  (V-v_3)(V-v_3^\sigma)(V-v_3^\tau)(V-v_3^{\sigma \tau}),\\
		\Phi_4(V) &=  (V-v_4)(V-v_4^\sigma)(V-v_4^\tau)(V-v_4^{\sigma \tau})\\
		& \quad \cdot (V-\overline{v_4})(V-\overline{v_4}^\sigma)(V-\overline{v_4}^\tau)(V-\overline{v_4}^{\sigma \tau}),
	\end{aligned}$$
	where $\overline{z}$ gives the conjugate of any complex number $z$, and the maps $\sigma$ and  $\tau$ change respectively the signs of $\sqrt{3}$, $\sqrt{5}$, and
	$v_1, v_2, v_3$, and $v_4$ are as follows:
	$$\begin{aligned}
		v_1 & =  \left( 3660\sqrt {5}-8190 \right) \sqrt {3}-6344\sqrt {5}+14175,\\
		v_2 & = \left( -420\sqrt {5}-990 \right) \sqrt {3}-784\sqrt {5}-1665,\\
		v_3 & = 315-440\I+ \left( 140-198\I \right) \sqrt {5},\\
		v_4 & =  \frac{\left( 3510-3300\I- \left( 1560-1485\I \right) \sqrt {5}
			\right) \sqrt {3} }{2} 
		+{\frac{6075}{2}}-2860\I \\
		&\quad -\left( 1350-1287\I 	\right) \sqrt {5}.
	\end{aligned}$$
	Hence, the $240$ roots of the  fundamental polynomial $\Phi(u^{12})$ of  $\Ee'_5$  are of the form
	$u=\zeta_{12}^\ell  v^{1/12} $	for $\ell=0, 1, \ldots,  11$, where
	$v$ varies on the set of $20$ roots of $\Phi(V)$.
	This means that the splitting field $\KK'_5$ of $\Ee'_5$ contains the field $k_0 (v_i^{1/12}: \ i=1,\cdots,4 )$.
	For each root of $\Phi(u)$, using the equations~\ref{eq4a}, one can determine the coefficients $a, b, c, d, e$  and hence a rational point
	$Q \in \Ee'_5(\KK'_5 (s))$ such that $sp_\infty (Q)=u$, 
	where $sp_\infty$ is the specializing  map  of $\Ee'_5(\KK'_5 (s))$ to the additive group  of $\KK'_5$.
	Indeed, it maps
	$Q \in \Ee'_5(\KK'_5(s))$ to the intersection point of the section $(Q)$ and the fiber over $\infty$ which lies in the smooth part of the additive singular fiber $\pi^{-1} (\infty)$.
	
	Since    $\Ee'_5$ has no reducible fiber  and all the $240$ sections $Q_j$'s  corresponding to the roots of  $\Phi(u)$
	are   points in $\Ee'_5(\KK'_5(s))$ with polynomial coordinates, we have
	$\langle Q_j, Q_j \rangle = 2$   and
	$\langle Q_{j_1}, Q_{j_2}\rangle
	= 1- ( Q_{j_1} \cdot Q_{j_2})$, where  $(Q_{j_1} \cdot Q_{j_2})$ denotes the intersection number
	for any $1\leq j_1 \neq j_2 \leq 2$.
	Assuming  $x_j=x(Q_j)$ and $y_j=y(Q_j)$,  the number
	$(Q_{j_1} \cdot Q_{j_2})$ can be computed by the following formula:
	\begin{align*}
		%\label{eqin}	
		(Q_{j_1} \cdot Q_{j_2})
		& =\deg(\gcd(x_{j_1}- x_{j_2}, y_{j_1}- y_{j_2})) \\
		& + 	\min \{2- \deg(x_{j_1}- x_{j_2}), 3- \deg(y_{j_1}- y_{j_2})\}.
	\end{align*}
	
	Using this formula and  determinant condition, we  obtain a subset of  eight points   with unimodular height paring matrix.
	In order to describe those points,  we consider the followings roots of $\Phi (V)$:
	%we let $\epsilon_4, \epsilon_5, \epsilon_6$  be as in
	% Table \eqref{tab:notation}.
	% the statement of Theorem~\ref{main-c1}.
	%In the next subsection,  we determine exactly the splitting field $\KK'_5$ of $\Ee'_5$  over $\Q(s)$ and then write down the coordinates of  eight generators of
	%$ \Ee'_5(\KK'_5(s))$.
	%\subsection{The splitting field $\KK'_5$ and eight generators of $\Ee'_5(\KK'_5)$}
	%\label{indps}
	%The Mordell--Weil lattice  $\Ee'(\CC(s))$ which is known to be isomorphic to the root lattice $E_8$, see for example \cite{Chahal2000, Usui2008}.
	%	Then, one can check that
	$$	\begin{aligned}[b]
		v_1& = \sqrt {5} \zeta_{12}^6 \epsilon_4^6 \epsilon_5^{-6}, &
		v_1^\sigma &=\sqrt {5} \zeta_{12}^6 \epsilon_4^{-6} \epsilon_5^6, &
		v_1^\tau &=\sqrt {5}\zeta_{12}^6  \epsilon_4^6 \epsilon_5^6, \\
		v_1^{\sigma \tau} &=\sqrt {5} \zeta_{12}^6  \epsilon_4^{-6} \epsilon_5^{-6}, &
		v_2  &= \sqrt {5} \zeta_{12}^6  \epsilon_4^6 \epsilon_5^{-6}  \epsilon_6^{12},  &
		v_2^\sigma &=\sqrt {5} \zeta_{12}^6 \epsilon_4^{-6} \epsilon_5^6   \epsilon_6^{12}, \\
		v_2^\tau &=\sqrt {5} \zeta_{12}^6 \epsilon_4^6 \epsilon_5^6 \epsilon_6^{-12},  &
		v_2^{\sigma \tau} &=\sqrt {5} \zeta_{12}^6  \epsilon_4^{-6}  \epsilon_5^{-6} \epsilon_6^{-12}.
	\end{aligned}$$
	Hence,  the roots $v_1, v_1^\sigma, v_1^\tau,  v_1^{\sigma \tau}, v_2, v_2^\sigma, v_2^\tau$ and $v_2^{\sigma \tau}$
	correspond  to  the following eight points:
	\begin{equation*}
		%\label{points}
		Q_j=\left(\frac{s^2 +a_j s +b_j}{u_j^2}, \frac{s^3 + c_j s^2 + d_j s + e_j}{u_j^3} \right) \in \Ee'_5(\KK'_5(s)),
	\end{equation*}
	for $j=1, \ldots, 8$, where  $a_j, b_j, c_j, d_j, e_j$ are  given in \cite[check-5]{k3-codes} and $u_j$'s are  given as follows:
	$$	\begin{aligned}[b]
		u_1 & = v_1^{\frac{1}{12}}, &
		%= \I \, \sqrt[24]{5} \, \left( \epsilon_4 \epsilon_6^{-1}\right)^{\frac 12},&
		u_2  &= (v_1^\sigma)^{\frac{1}{12}}, &
		%= \I \, \sqrt[24]{5}\, \left( \epsilon_4^{-1} \epsilon_6 \right)^{\frac 12},\\
		u_3  &= (v_1^\tau)^{\frac{1}{12}}, &
		%= \I \, \sqrt[24]{5} \,  \left(\epsilon_4 \epsilon_6\right)^{\frac 12},&
		u_4  &= (v_1^{\sigma \tau})^{\frac{1}{12}}, \\
		% = \I \, \sqrt[24]{5} \, \left( \epsilon_4 \epsilon_6\right)^{-\frac 12},\\
		u_5 & = v_2^{\frac{1}{12}}, &
		%= \I \, \sqrt[24]{5} \, \left( \epsilon_4\epsilon_5^2 \epsilon_6^{-1}\right)^{\frac 12},  &
		u_6  &= (v_2^\sigma)^{\frac{1}{12}}, &
		%= \I \, \sqrt[24]{5} \, \left( \epsilon_4^{-1}\epsilon_5^2 \epsilon_6\right)^{\frac 12}, \\
		u_7  &= (v_2^\tau)^{\frac{1}{12}}, &
		%	= \I \, \sqrt[24]{5}\, (\epsilon_4\epsilon_5^{-2} \epsilon_6)^{\frac 12}, &
		u_8  &= (v_2^{\sigma \tau})^{\frac{1}{12}}.
		%= \I \, \sqrt[24]{5} \,\left( \epsilon_4 \epsilon_5 \epsilon_6\right)^{-\frac 12}.
	\end{aligned}$$
	Applying the specialization map $sp_\infty: \Ee'_5(\KK'_5) \rightarrow (\KK'_5)^+$ to these points and dividing the images by $u_1$,
	we obtain the following subset of the  field  $\KK'_5$,
	$$\left\{1, \frac{u_j}{u_1}: j= 2, \ldots,  8\right\}=\left\{ 1, \epsilon_4^{-1} \epsilon_6,
	\epsilon_6,  \epsilon_4^{-1},  \epsilon_5,  \epsilon_4^{-1}\epsilon_5 \epsilon_6,  \epsilon_5^{-1} \epsilon_6,  \epsilon_4^{-1} \epsilon_5^{-1} \right\},$$
	which is easy to see that they are linearly independent over $\Q$. Thus, the  points
	$Q_1, \ldots ,  Q_8$  form a linearly independent subset   generating a sublattice of rank $8$ in $\Ee'_5(\KK'_5(s)).$
	%	and have an unimodular Gram matrix  as follows:
	The Gram matrix  of these eight points  is equal to the following unimodular matrix:
	{\small 	$$R'_5=\begin{pmatrix}
			2 & 0 & 0 & -1 & 0 & 0 & 0 & 0 \\ \noalign{\medskip}
			0 & 2 & -1 & 0 & 0 & 0 & 1 & -1  \\ \noalign{\medskip}
			0 & -1 & 2 & 0 & 0 & 0 & 0 & 0 \\ \noalign{\medskip}
			-1 & 0 & 0 & 2 & -1 & -1 & 0 & 0 \\ \noalign{\medskip}
			0 & 0 & 0 & -1 & 2 & 0 & 0 & -1 \\ \noalign{\medskip}
			0 & 0 & 0 & -1 & 0 & 2 & -1 & 0 \\ \noalign{\medskip}
			0 & 1 & 0 & 0 & 0 & -1 & 2 & 0 \\ \noalign{\medskip}
			0 & -1 & 0 & 0 & -1 & 0 & 0 & 2
		\end{pmatrix}. $$}
	%
	%	By a direct searching between $20$ roots of the polynomial $\Phi(V)$, we find
	%	$89$ other  $8$-tuples with unimodular Gram matrix.
	Thus, 	the points  $Q_j$'s for  $ \ j= 1, \ldots, 8$ generate the whole group  $\Ee'_5(\KK'_5(s))$ as desired.
	%	Hence,  the specializing map is an isomorphism and the splitting
	%	field $\KK'_5$ is obtained by adjoining one of the $u_j$'s, for example $u_3$, to the field $k_0$. 
	%	In other words, we have $\KK'_5=k_0(u_3)$ which contains $\Q\left(\zeta_{12}, 5^{\frac{1}{24}}, \left(\epsilon_4 \epsilon_5\right)^{\frac 1 2}\right).$	
	Using Pari/GP, we   obtained the minimal defining polynomial of $\KK'_5$ is a number  field 
	defined by a polynomial $g'_5(x)$ of degree 96  given in \cite[min-pols]{k3-codes}.
	\iffalse
	as follows:
	\begin{align}
		g'_5(x)&= x^{96} + 110880x^{84} + 2930852960x^{72} - 7778609568000x^{60} \notag \\
		&\quad + 12574087596678400x^{48} - 9644688087731712000x^{36} \notag \\
		&\quad + 4571547377069445120000x^{24}   - 17794702679917363200000x^{12} \notag \\
		&\quad + 51322966254895104000000.
		\label{g'5}
	\end{align}
	\fi
	Therefore,   the proof of Theorem~\ref{main-c1} is  finished.
	\subsection{Proof of Theorem~\ref{main-c}}
	%++++++++++++++++++++++++++++++++++++++++=++++
	
	First, we note that the splitting field of the elliptic $K3$ surface $\Ee_5$ over $\Q(t)$
	is equal to $\KK_5=\KK'_5(\zeta_5)$ where $\KK'_5$ is the splitting field of  $\Ee'_5$  over $\Q(s)$. The number field  $\KK_5$ is defined by a polynomial $g_5(x)$ of degree 192 with huge coefficients given in \cite{k3-codes}. Indeed the  $\KK_5$ is the compositum of the polynomial $g'_5(x)$  and
	the  cyclotomic polynomial of order 5, since
	$x^5-1=(x-1) (x^4 + x^3 + x^2 + x + 1)$.

	Letting $s= t+ 1/t$, the rational elliptic surface
	$\Ee'_5$ over $\KK_5 (s)$ is isomorphic to
	$\Ee_5$ over  $\KK_5 (t)$ as a quadratic extension of $\KK_5 (s)$. Hence, the independent generators
	$$ Q_j=\left(\frac{s^2 + a_j s + b_j}{u_j^2},\frac{s^3+ c_j s^2 + d_j s + e_j}{u_j^3}  \right) $$
	of $\Ee'(\KK_5 (s))$ leads  to the   points  $P_j=(x_j(t), y_j(t))\in \Ee_5(\KK_5 (t))$ 
	of the form	given in the statement of  Theorem~\ref{main-c} with
	%
	%$$\begin{aligned}
		%		x_j(t) &= \frac{t^4 + a_j t^3 + (b_j +2) t^2 + a_j t +1}{u_j^2  t^2}, \\
		%		y_j(t) &=\frac{t^6 +c_j t^5+ (d_j +3) t^4 +(2 c_j +e_j) t^3 +(d_j +3) t^2 + c_j t +1}{u_j^3  t^3},
		%	\end{aligned}$$
	the constants $a_j, b_j, c_j, d_j, e_j$  and  $u_j$'s for $j=1,\ldots, 8,$ provided in 
	\cite{k3-codes}.
	%	 Subsection~\ref{fpe0}.	
	
	Furthermore,
	%by fixing $\zeta_5 =\left( \sqrt{5}-1+\mathrm{i} \sqrt{2}\, \sqrt{5+\sqrt{5}}\right)/4$,
	by letting $s= \zeta_5 t+ \frac{1}{\zeta_5 t}$ and the same argument as above, 
	we obtain points $P_{j+8} =(x_{j}(\zeta_5 t ), y_{j}(\zeta_5 t ))$ 
	%	with coordinates
	%	\begin{align*}
		%		x_j(t) &= \frac{\zeta_5^4\, t^4 + a_j\, \zeta_5^3\, t^3 + (b_j +2) \zeta_5^2\, t^2 + a_j \zeta_5\, t +1}{u_j^2 \, \zeta_5^2\, t^2}, \\
		%		y_j(t) &=\frac{\zeta_5\, t^6 +c_j t^5+ (d_j +3)\, \zeta_5^4\, t^4 +(2 c_j +e_j)\, \zeta_5^3\, t^3 +(d_j +3)\,\zeta_5^2\, t^2 + c_j\, \zeta_5\, t +1}{u_j^3 \,\zeta_5^3\, t^3},
		%	\end{align*}
	for $j=1,\ldots, 8$. We note that the points $P'_j= \left(t^2 x(P_j), t^3  y(P_j) \right)$ belong to
	the Mordell--Weil lattice of the elliptic $K3$ surface $\Ee: y^2= x^3+ t (t^{10}+1)$, which is birational to $\Ee_5$	over $\KK_5(t)$.
	%	We refer the reader to   \cite{Usui2008} for more details.
	Since   $P'_j$'s have no intersection with the zero section of $\Ee$, we have
	$\langle P'_j, P'_j \rangle = 4$, and 
	$\langle P'_{j_1}, P'_{j_2}\rangle = 2- (P'_{j_1} \cdot P'_{j_2})$,
	for  $1\leq j_1 \neq j_2 \leq 16$. The intersection number
	$(P'_{j_1} \cdot P'_{j_2})$ can be computed by
	\begin{align*}
		( P'_{j_1} \cdot P'_{j_2}) &= \deg\left(\gcd (x_{j_1}- x_{j_2} , y_{j_1}- y_{j_2})\right) \notag \\
		&\quad + \min \left\{4- \deg( x_{j_1}- x_{j_2}) , 6- \deg( y_{j_1}- y_{j_2})\right\}
	\end{align*}
	Thus,  we obtain the following Gram matrix $R_5$ of the height pairing  for points $P'_j$'s and hence $P_j$'s:
	%	
	%$$R_5:=\left(
	%\begin{array}{cccccccccccccccc}
	%	4 & 0 & 0 & 0 & 0 & -2 & 0 & 2 & 0 & 0 & 0 & 0 & -1 & 1 & -1 & -1
	%	\\ \noalign{\medskip}
	%	0 & 4 & -2 & 2 & -2 & 0 & 0 & 0 & 0 & 0 & 0 & -1 & -1 & 0 & -1 & 1
	%	\\ \noalign{\medskip}
	%	0 & -2 & 4 & -2 & 0 & -2 & 0 & 0 & 0 & 0 & 0 & 0 & 0 & 0 & 0 & -1
	%	\\ \noalign{\medskip}
	%	0 & 2 & -2 & 4 & 0 & 0 & 2 & 0 & 0 & -1 & 0 & -2 & 0 & 1 & -2 & 1
	%	\\ \noalign{\medskip}
	%	0 & -2 & 0 & 0 & 4 & 0 & 2 & 0 & -1 & -1 & 0 & 0 & 0 & 1 & 0 & -1
	%	\\ \noalign{\medskip}
	%	-2 & 0 & -2 & 0 & 0 & 4 & -2 & -2 & 1 & 0 & 0 & 1 & 1 & -2 & 2 & 2
	%	\\ \noalign{\medskip}
	%	0 & 0 & 0 & 2 & 2 & -2 & 4 & 0 & -1 & -1 & 0 & -2 & 0 & 2 & -2 & -1
	%	\\ \noalign{\medskip}
	%	2 & 0 & 0 & 0 & 0 & -2 & 0 & 4 & -1 & 1 & -1 & 1 & -1 & 2 & -1 & -2
	%	\\ \noalign{\medskip}
	%	0 & 0 & 0 & 0 & -1 & 1 & -1 & -1 & 4 & 0 & 0 & 0 & 0 & -2 & 0 & 2
	%	\\ \noalign{\medskip}
	%	0 & 0 & 0 & -1 & -1 & 0 & -1 & 1 & 0 & 4 & -2 & 2 & -2 & 0 & 0 & 0
	%	\\ \noalign{\medskip}
	%	0 & 0 & 0 & 0 & 0 & 0 & 0 & -1 & 0 & -2 & 4 & -2 & 0 & -2 & 0 & 0
	%	\\ \noalign{\medskip}
	%	0 & -1 & 0 & -2 & 0 & 1 & -2 & 1 & 0 & 2 & -2 & 4 & 0 & 0 & 2 & 0
	%	\\ \noalign{\medskip}
	%	-1 & -1 & 0 & 0 & 0 & 1 & 0 & -1 & 0 & -2 & 0 & 0 & 4 & 0 & 2 & 0
	%	\\ \noalign{\medskip}
	%	1 & 0 & 0 & 1 & 1 & -2 & 2 & 2 & -2 & 0 & -2 & 0 & 0 & 4 & -2 & -2
	%	\\ \noalign{\medskip}
	%	-1 & -1 & 0 & -2 & 0 & 2 & -2 & -1 & 0 & 0 & 0 & 2 & 2 & -2 & 4 & 0
	%	\\ \noalign{\medskip}
	%	-1 & 1 & -1 & 1 & -1 & 2 & -1 & -2 & 2 & 0 & 0 & 0 & 0 & -2 & 0 & 4
	%\end{array}
	%\right).	$$
	
	{\small  
		$$R_5=
		\left( \begin {array}{cccccccccccccccc}
		4 & 0 & 0 & 2 & 0 & 0 & 0 & 0 & -2 & 0 & 0 & 0 & 2 & 0 & 0 & 1 	\\ \noalign{\medskip}
		0 & 4 & 2 & 0 & 0 & 0 & -2 & 2 & 0 & -2 & -1 & 1 & 0 & 2 & 2 & -1 	\\ \noalign{\medskip}
		0 & 2 & 4 & 0 & 0 & 0 & 0 & 0 & 0 & -1 & 0 & 0 & 0 & 1 & 2 & 0 	\\ \noalign{\medskip}
		2 & 0 & 0 & 4 & 2 & 2 & 0 & 0 & 0 & 1 & 0 & 0 & 1 & 0 & 0 & 2 	\\ \noalign{\medskip}
		0 & 0 & 0 & 2 & 4 & 0 & 0 & 2 & 2 & 0 & 0 & 1 & 0 & 0 & 0 & 0 	\\ \noalign{\medskip}
		0 & 0 & 0 & 2 & 0 & 4 & 2 & 0 & 0 & 2 & 1 & 0 & 0 & 0 & -1 & 1 	\\ \noalign{\medskip}
		0 & -2 & 0 & 0 & 0 & 2 & 4 & 0 & 0 & 2 & 2 & 0 & 0 & -1 & -2 & 0 	\\ \noalign{\medskip}
		0 & 2 & 0 & 0 & 2 & 0 & 0 & 4 & 1 & -1 & 0 & 2 & 0 & 1 & 0 & -2 	\\ \noalign{\medskip}
		-2 & 0 & 0 & 0 & 2 & 0 & 0 & 1 & 4 & 0 & 0 & 2 & 0 & 0 & 0 & 0 	\\ \noalign{\medskip}
		0 & -2 & -1 & 1 & 0 & 2 & 2 & -1 & 0 & 4 & 2 & 0 & 0 & 0 & -2 & 2 	\\ \noalign{\medskip}
		0 & -1 & 0 & 0 & 0 & 1 & 2 & 0 & 0 & 2 & 4 & 0 & 0 & 0 & 0 & 0 	\\ \noalign{\medskip}
		0 & 1 & 0 & 0 & 1 & 0 & 0 & 2 & 2 & 0 & 0 & 4 & 2 & 2 & 0 & 0 	\\ \noalign{\medskip}
		2 & 0 & 0 & 1 & 0 & 0 & 0 & 0 & 0 & 0 & 0 & 2 & 4 & 0 & 0 & 2 	\\ \noalign{\medskip}
		0 & 2 & 1 & 0 & 0 & 0 & -1 & 1 & 0 & 0 & 0 & 2 & 0 & 4 & 2 & 0 	\\ \noalign{\medskip}
		0 & 2 & 2 & 0 & 0 & -1 & -2 & 0 & 0 & -2 & 0 & 0 & 0 & 2 & 4 & 0 	\\ \noalign{\medskip}
		1 & -1 & 0 & 2 & 0 & 1 & 0 & -2 & 0 & 2 & 0 & 0 & 2 & 0 & 0 & 4
	\end{array} \right).$$ 
}

One can check that its determinant is equal to $5^4$ as desired, which shows that the   points
$P_j$'s for $j=1, \ldots, 16$ form a set of independent generators of $\Ee_5 $	over $\KK_5(t)$.
We refer the reader to see \cite[check-5]{k3-codes} for the computations of this section, and 
\cite[Points-5]{k3-codes} the list of 16 points in $\Ee_5  (\KK_5(t))$.

%=======================================================
\section{The case of $\Ee_6$}
%=================================================
\label{case6}

We prove Theorem~\ref{main-d} on the  elliptic $K3$ surface
$\Ee_6: y^2= x^3 + t^6 + 1/t^6$ over $\CC(t)$ 	in this section.
To do this, first we  determine  the splitting field  $\KK'_6$ and find a set of independent generators for
the Mordell--Weil lattice of the rational elliptic surface
$\Ee'_6: y^2= x^3+ f_6(s)$ where
$$f_6(s)= s^6- 6 s^4 + 9 s^2-2= (s^2-2)(s^4- 4 s^2 +1).$$
To simplify the computations, we set $\ts= s-\sqrt{2}$   to obtain  the  rational elliptic surface
\begin{equation}
	\label{eq6-0}
	\widetilde{\Ee}'_6: y^2= x^3+ f'_6(\ts), 
\end{equation}
where 
$$ f'_6(\ts)=  \ts (\ts- 2\sqrt{2})(\ts^2 - \sqrt{2} \ts- 1)(\ts^2 - 3 \sqrt{2} \ts + 3).$$
It is easy to see that  $\widetilde{\Ee}'_6$ is  birational to  $\Ee'_6$ over $\Q ( \sqrt{2})$.
Since $ \widetilde{\Ee}'_6 (\CC(\ts)) \iso  \Ee'_6 (\CC(s)) \iso E_8$, there exist exactly  $240$ points
in $\widetilde{\Ee}'_6 (\CC(\ts))$
of the form
\begin{equation}
	\label{eq6-01}
	\tilde{Q}=\left( a \ts^2 +b \ts + g,\  c \ts^3 + d \ts^2 + e \ts + h  \right)
\end{equation}
corresponding to the   points  $Q=(x,y) \in \Ee'_6 (\CC(s))$ with
\begin{equation}
	\label{eq6-1}
	\begin{aligned}[b]
		x(s)&= a s^2 +(b - 2 \sqrt{2} a) s + g +(2 a - \sqrt{2} b),  \\ 
		y(s)&= c s^3 + (d-3  \sqrt{2} c) s^2 + (6 c - 2 \sqrt{2}d + e) s
		+ h  -\sqrt{2}(2  c -  \sqrt{2} d + e).
	\end{aligned}
\end{equation}

It is clear that the splitting  field $\KK'_6$ of $\Ee'$ is a quadratic extension by $\sqrt{2}$ of the splitting field of   $\widetilde{\Ee}'_6$, which is denoted by
$\widetilde{\KK}'_6$ and contains $\Q(\sqrt{2})$ as  a subfield.
We have the following result on  $\widetilde{\KK}'_6$ and the set of independent generators of
$\widetilde{\Ee}'_6 (\widetilde{\KK}'_6(\ts))$.

\begin{thm}
	\label{main-d1}
	The splitting field  $\widetilde{\KK}'_6$ of 
	$ \widetilde{\Ee}'_6 $  is a number field  of degree 96, with a minimal defining polynomial  given in \cite{k3-codes},	containing 
	$\Q(\mathrm{i}, \beta_0, \beta_1, {u_7}^{\frac{1}{2}})$,
	where   $u_8$ is given by \eqref{6eqn2} below, $\beta_0$ and $\beta_1$  as in Table~\ref{tab:notation}.
	%	$\beta_0:= 2^{\frac{1}{6}},
	%	\beta_1:= (3+ 2 \sqrt{3})^{\frac{1}{4}}.$
	
	Moreover, the   lattice  $\widetilde{\Ee}'_6 (\widetilde{\KK}'_6(\ts))$ is generated by
	\begin{align*}
		\tilde{Q}_j&=\left( a_j \ts^2 +b_j \ts + u_j^2,\  c_j \ts^3 + d_j \ts^2 + e_j \ts + u_j^3  \right), \ (j=1,\ldots, 8)
	\end{align*}
	where  $a_j, b_j, c_j, d_j, e_j$ and $u_j$ are given in  \cite[check-6]{k3-codes}.
\end{thm}

In the next   subsection, we prove the above theorem and in Subsection~\ref{last} we provide the complete proof of 	Theorem~\ref{main-d}.

\subsection{Proof of Theorem~\ref{main-d1}}
\label{coefs6}
We first determine the fundamental polynomial of $\widetilde{\Ee}'_6$  over $\Q(\sqrt{2})$.
Substituting the coordinates of $\widetilde{Q} \in \widetilde{\Ee}'_6 (\widetilde{\KK}'_6(\ts))$, given by \eqref{eq6-01}, in the equation~\eqref{eq6-0} of $\widetilde{\Ee}'_6$ and letting $g = u^2, h = u^3$, we get  the following six relations:
\[
\begin{aligned} 
	\label{6eqn1}
	a^3-c^2+1&=0,  &
	3a^2b-2cd-6\sqrt{2} &=0, \\
	3a^2u^2+3ab^2-2ce-d^2+24&=0,&
	6abu^2-2cu^3+b^3-2de-16\sqrt{2}&=0, \\
	3au^4+3b^2u^2-2du^3-e^2-3&=0, &
	3bu^4-2eu^3+6\sqrt{2}&= 0.
\end{aligned}
\]
Using Maple , one can compute   the
fundamental polynomial  $\Phi (u)$ of the ideal generated by  the above equations,
which is a polynomial of degree $240$ in terms of $u$ up to a constant.
By taking $v=u^2$, we obtain a polynomial $\Phi (v)$  in $\Z[v]$ which can be decomposed into nine irreducible factors, namely,	
$$\Phi (v) = \prod_{i=1}^{9}\Phi_i(v).$$
The first six factors of $\Phi(v)$ can be decomposed 
as follows:
$$	\begin{aligned}
	\Phi_1(v) &= v^3-2 = (v - \beta_0^2 )(v - \beta_0^2  \zeta_3)(v - \beta_0^2  \zeta_3^2),  \\
	\Phi_2(v) &= v^4-6 v^2-3= (v-\beta_1^2)(v+\beta_1^2)(v-\I  \beta_2^2)(v+ \I \beta_2^2),   \\
	\Phi_3(v) &= v^3+12 v^2+12 v+6= (v- v_{31})(v- v_{32})(v- v_{33}), \\
	\Phi_4(v) &= v^8+6 v^6+39 v^4-18 v^2+9 \\
	& = (v+ \zeta_3 \beta_1^2) (v- \zeta_3 \beta_1^2) (v+ \zeta_6 \beta_1^2)(v- \zeta_6 \beta_1^2) \\
	&  \quad \times 
	(v+ \zeta_{12} \epsilon_1^{'} \beta_1^2)(v- \zeta_{12} \epsilon_1^{'} \beta_1^2)(v+\zeta_{12}^{11} \epsilon_1^{'} \beta_1^2)(v-\zeta_{12}^{11} \epsilon_1^{'} \beta_1^2),
\end{aligned}$$	
$$	\begin{aligned}
	\Phi_5(v) &= v^6-12 v^5+132 v^4-132 v^3+72 v^2-72 v+36,\\
	& =(v- v_{51})(v- v_{51}^\gamma)(v- v_{52})(v- v_{52}^\gamma)(v- v_{53})(v- v_{53}^\gamma),\\
	\Phi_6(u) &= v^8-48 v^7+168 v^6-912 v^5+1272 v^4-1152 v^3+864 v^2-576 v+144\\
	& =(v-v_{61})(v-  v_{61}^\gamma)(v-v_{62})(v-v_{62}^\gamma)(v-v_{63}) (v-v_{63}^\gamma)(v-v_{64})(v-v_{64}^\gamma), 		
\end{aligned}$$
where 
$$
\begin{aligned}
	v_{31} &= -( 2 \beta_0^4 + 3 \beta_0^2  + 4), \
	v_{32}=  2 \beta_0^4 \zeta_{6}^{5} + 3 \beta_0^2  \zeta_{6}  - 4, \\
	v_{33}& =  2 \beta_0^4 \zeta_{6}  + 3 \beta_0^2  \zeta_{6}^{5} - 4, \\
	v_{51} &= (2 \beta_0^4 + 3 \beta_0^2 +4) \zeta_{6},\
	v_{52}=2 \beta_0^4  \zeta_{6} - 3 \beta_0^2 +4  \zeta_{6}^{5},\\
	v_{53}& =- 2 \beta_0^4  + 3 \beta_0^2  \zeta_{6} +4  \zeta_{6}^{5}, \\			
	v_{61} & = \sqrt{3} \beta_1^{4}+2 \mathrm{i} \sqrt{2} \beta_1^{3}-(2 \sqrt{3}+1) \beta_1^{2}
	-\mathrm{i} \sqrt{2} (\sqrt{3}+3) \beta_1, \\
	v_{62}& = \sqrt{3} \beta_1^{4}+2 \sqrt{2} \beta_1^{3}+(2 \sqrt{3}+1) \beta_1^{2}+\sqrt{2} (\sqrt{3}+3) \beta_1, \\
	v_{63}& =  -\sqrt{3} \beta_1^{4}+\left(\mathrm{i} +1\right) \sqrt{2} (3 \sqrt{3}-5) \beta_1^{3}+\mathrm{i} (5 \sqrt{3}-8) \beta_1^{2}\\
	&+(\mathrm{i}-1) \sqrt{2} (-3+2 \sqrt{3}) \beta_1 +12,\\
	v_{64}& = -\sqrt{3} \beta_1^{4}+(\mathrm{i}-1) \sqrt{2} (3 \sqrt{3}-5) \beta_1^{3}-\mathrm{i} (5 \sqrt{3}-8) \beta_1^{2}\\
	& +(\mathrm{i}+1) \sqrt{2} (-3+2 \sqrt{3}) \beta_1 +12,
\end{aligned}$$
and $\gamma$  changes the sign of $\sqrt{2}$.
One can check that the other three factors of $\Phi(v)$, say $\Phi_7(v)$,  $\Phi_8(v)$  and $\Phi_9(v)$  of  degrees 16,  24 and  48 respectively,   can be completely decomposed over  $\Q(\mathrm{i}, \beta_0, \beta_1)$. 	 
For example,  the seventh factor  is    following degree $16$ polynomial:
$$ 	\begin{aligned}
	\Phi_7(v) &=     v^{16}+48 v^{15}+2136 v^{14}+6240 v^{13}-16824 v^{12}+32256 v^{11} +564480 v^{10}\\
	& \quad 	+815040 v^9+477360 v^8 -6912 v^7-248832 v^6-338688 v^5\\
	&\quad -100224 v^4+165888 v^3+207360 v^2+82944 v+20736,
\end{aligned}$$
and one of its roots is equal to
\begin{equation}
	\label{6eqn2}
	v_{7}=\frac{1}{2} \left(v_{73} \beta_1^3 + v_{72} \beta_1^2 +v_{71} \beta_1+ v_{70}\right)
\end{equation}
where
%\begin{equation}
%		\label{6eqn2}
%\begin{align*}
%%	v_{7}& =\frac{1}{2} \left(v_{73} \beta_1^3 + v_{72} \beta_1^2 +v_{71} \beta_1+ v_{70}\right) & & \\
%	v_{70} & = 3 \big((1+ 2\mathrm{i})\sqrt{3}  -  (2+ 3\mathrm{i} ) \big), &
%	v_{71} & = - \beta_0^3 \big((5\mathrm{i}-1)\sqrt{3}  +  (3-9\mathrm{i} )\big),\\
%	v_{72} & =  (  5\mathrm{i} -8)\sqrt{3} +  (15- 8\mathrm{i} ),&
%	v_{73} & = 2 \beta_0^3 \big((4+ \mathrm{i})\sqrt{3}  -  (7+ 2\mathrm{i} )\big).
%\end{align*}
$$\begin{aligned}[b]
	v_{70} &= 3 \big((1+ 2\mathrm{i})\sqrt{3}  -  (2+ 3\mathrm{i} ) \big), &
	v_{71}  &= - \beta_0^3 \big((5\mathrm{i}-1)\sqrt{3}  +  (3-9\mathrm{i} )\big), \\ 
	v_{72}  &=  (  5\mathrm{i} -8)\sqrt{3} +  (15- 8\mathrm{i} ), &
	v_{73}  &= 2 \beta_0^3 \big((4+ \mathrm{i})\sqrt{3}  -  (7+ 2\mathrm{i} )\big).
\end{aligned}$$
We cite  \cite[check-6]{k3-codes}   to   see    the complete decomposition of all factors of $\Phi(v)$.
Thus, the field $\KK_6^{'}$  is an extension of  $\Q(\mathrm{i}, \beta_0, \beta_1)$, 
with a defining minimal polynomial of degree 96 with huge coefficients, see \cite[min-pols]{k3-codes}.

By   a direct searching for eight roots  between 240 root of  	$\Phi(u)$  determinant conditions on the
Gram matrix of  corresponding points,  we find
the following roots:
$$
\begin{aligned}[b]
	u_1 & =\beta_0, &
	u_2 & = \zeta_6 \beta_0, &
	u_3 & =   \beta_1,&
	u_4 & =  \zeta_8 \beta_1,\\
	u_5 & =   v_{32}^{\frac{1}{2}},&
	u_6 & =  \zeta_{12} \beta_1,&
	u_7 & =  v_{61}^{\frac{1}{2}}, &
	u_8 & =  v_{7}^{\frac{1}{2}}.
\end{aligned}
$$
These  provide  the eight points  generating  
$\widetilde{\Ee}'_6 (\widetilde{\KK}'_6(\ts))$ as follows:
\begin{align*}
	\tilde{Q}_j&=\left( a_j \ts^2 +b_j \ts + g_j,\  c_j \ts^3 + d_j \ts^2 + e_j \ts + h_j  \right), \ (j=1,\ldots, 8)
\end{align*}
where  $a_j, b_j, c_j, d_j, g_i, h_i$ are given in \cite[check-6]{k3-codes}.

Applying the specialization map $sp_0: \Ee'_6(\KK'_6(\ts)) \rightarrow (\KK'_6)^+$,
defined by
$$P \mapsto sp_0(P)= \frac{1}{u}= \frac{x(P)}{y(P)}\bigg|_{\ts=0},$$
to the  above points  and multiplying  the images by $u_1$,
we obtain
$$\left\{1, \frac{u_j}{u_1}: j= 2, \ldots,  8\right\} \subset \KK'_6,$$
which can be checked that they are linearly independent over $\Q$. Thus the  points
$\tilde{Q}_1, \ldots ,  \tilde{Q}_8$  form a linearly independent subset   generating a sublattice of rank $8$ in $\Ee'_6(\KK'_6(\ts)).$
Moreover, the Gram matrix  of the  points
$\tilde{Q}_1, \ldots, \tilde{Q}_8$
is equal to the following unimodular matrix:
$$R'_6=%\left(\begin {array}{rrrrrrrr}
\begin{pmatrix}
	2&1&0&0&1&0&0&0\\ \noalign{\medskip}
	1&2&0&0&1&0&0&0\\ \noalign{\medskip}
	0&0&2&0&0&-1&0&1\\ \noalign{\medskip}
	0&0&0&2&0&1&-1&1\\ \noalign{\medskip}
	1&1&0&0&2&0&0&1\\ \noalign{\medskip}
	0&0&-1&1&0&2&0&0\\ \noalign{\medskip}
	0&0&0&-1&0&0&2&0\\ \noalign{\medskip}
	0&0&1&1&1&0&0&2
\end{pmatrix}.
%\end {array} \right)
$$
Finally, one  can use   \eqref{eq6-1} to get
the   points  $Q_j=(x_j, y_j) \in \Ee'_6 (\KK'_6 (s))$ with
$$
\begin{aligned}[b]
	x_j(s)&= a_j s^2 +(b_j - 2 \sqrt{2} a_j) s + u_j^2 +(2 a_j - \sqrt{2} b_j),  \\
	y_j(s)&= c_j s^3 + (d_j-3  \sqrt{2} c_j) s^2 + (6 c_j - 2 \sqrt{2}d_j + e_j) s\\
	&  \quad	+ u_j^3  -\sqrt{2}(2  c_j -  \sqrt{2} d_j + e_j).
\end{aligned}
$$

\subsection{Proof of Theorem~\ref{main-d}}
\label{last}

The splitting field of the elliptic $K3$ surface $\Ee_6$ over $\Q(t)$
is equal to $\KK_6=\KK'_6(\zeta_{12}) = \KK'_6$,
where $\KK'_6$ is the splitting field of  $\Ee'_6: y^2=x^3+f_6 (s)$  over $\Q(s)$.

Letting $s= t+ 1/t$, the rational elliptic surface
$\Ee'_6$ over $\KK_6 (s)$ is isomorphic to
$\Ee_6$ over  $\KK_6 (t)$ as a quadratic extension of $\KK_6 (s)$. Hence, the eight independent generators
$ Q_j=(x_j(s), y_j(s))  \in \Ee'(\KK_6 (s)) $
give  the   points  $P_j=(x_j(t), y_j(t))\in \Ee_6(\KK_6 (t))$ with
\begin{align*}
	x_j(t)& = \frac{a_{j,4} t^{4}+a_{j,3} t^{3}+a_{j,2} t^{2}+ a_{j,1} t +a_{j,0}}{t^2}, \\
	y_j(t) &=\frac{b_{j,6}t^{6}+b_{j,5} t^{5}+b_{j,4} t^{4}+b_{j,3} t^{3}+b_{j,2} t^{2}+b_{j,1} t +b_{j,0}}{t^3},\\
	%\end{align*}
	%	
	%	\begin{align*}
		%a_{j,\,4}& =a_{j,\,0}=a_j, \ \ a_{j,\,3} = a_{j,\,1} =b_j-2 \sqrt{2}\,  a_j,  \ \
		%a_{j,\,2} =  g_j+4 a_j -\sqrt{2}\, b_j,
		%\\
		%b_{j,\,6}&=b_{j,\,6}=c_j, \ b_{j,\,5}=b_{j,\,1}= c_j +d_j -3 \sqrt{2}, \\
		%b_{j,\,4}&=b_{j,\,2}= d_j +9 c_j +e_j -2 \sqrt{2}, \
		%b_{j,\,3}=h_j-8 \sqrt{2}\, c_j -\sqrt{2}\, e_j +4 d_j.
	\end{align*}
	where
	$$\begin{aligned}
		a_{j,0}& =a_{j,4}=a_j, & a_{j,1}&= a_{j,3} =b_j-2 \sqrt{2}  a_j, \\
		a_{j,2} &=  u_j^2+4 a_j -\sqrt{2} b_j,&   &\\
		b_{j,0}&=b_{j,6}=c_j, & b_{j,1}&=b_{j,5}= c_j +d_j -3 \sqrt{2}, \\
		b_{j,2}&=b_{j,4}= d_j +9 c_j +e_j -2 \sqrt{2}, & 	b_{j,3}&=u_j^3-8 \sqrt{2} c_j -\sqrt{2} e_j +4 d_j.
	\end{aligned}$$
	The constants $a_j, b_j, c_j, d_j, e_j$  and  $u_j$'s for $j=1,\ldots, 8,$ are listed in the previous subsection.
	Furthermore,
	%by fixing $\zeta_5 =\left( \sqrt{5}-1+\mathrm{i} \sqrt{2}\, \sqrt{5+\sqrt{5}}\right)/4$,
	letting $s= \zeta_{12} t+ \frac{1}{\zeta_{12} t}$, same  as above, we obtain the
	points $P_{j+8} =(x_{j+8}(t), y_{j+8}(t))$ with coordinates
	$$\begin{aligned}[b]
		x_{j+8}(t)& = \frac{a_{{j+8},4} t^{4}+a_{{j+8},3} t^{3}+a_{{j+8},2} t^{2}+ a_{{j+8},1} t +a_{{j+8},0}}{\zeta_{12}^2  t^2}, \\
		y_{j+8}(t) &=\frac{b_{{j+8},6}t^{6}+b_{{j+8},5} t^{5}+b_{{j+8},4} t^{4}+b_{{j+8},3} t^{3}+b_{{j+8},2} t^{2}+b_{{j+8},1} t +b_{{j+8},0}}{\zeta_{12}^3  t^3},
	\end{aligned}$$
	where
	%	$$\begin{aligned}[b]
		%		a_{{j+8},4} & =\zeta_3 a_{{j+8},0} = \zeta_3  a_j, & 
		%		a_{{j+8},3} & =\zeta_6 a_{{j+8},1} = \zeta_6 2 \sqrt{2} a_j+b_j,\\
		%		a_{{j+8},2} & = \zeta_6 ( a_j + \sqrt{2} b_j + g_j), &
		%		b_{{j+8},6} & = -b_{{j+8},0} = c_j,\\
		%		b_{{j+8},5} & =b_{{j+8},1}=\zeta_{12}^5 (3 \sqrt{2} c_j + d_j), &
		%		b_{{j+8},4} & =b_{{j+8},2} =  \zeta_3 ( 9 c_j +2 \sqrt{2} d_j + e_j ),\\
		%		b_{{j+8},3} & = \mathrm{i} (8\sqrt{2} c_j + 4 d_j + \sqrt{2} e_j + h_j ),
		%	\end{aligned}$$
	$$\begin{aligned}[b]
		a_{{j+8},4} & =\zeta_3 a_{{j+8},0} = \zeta_3  a_j, & 
		a_{{j+8},3} & =\zeta_6 a_{{j+8},1} = \zeta_6 2 \sqrt{2} a_j+b_j,\\
		a_{{j+8},2} & = \zeta_6 ( a_j + \sqrt{2} b_j + g_j), & & \\
		b_{{j+8},6} & = -b_{{j+8},0} = c_j, &
		b_{{j+8},5} & = b_{{j+8},1}=\zeta_{12}^5 (3 \sqrt{2} c_j + d_j), \\
		b_{{j+8},4} & =b_{{j+8},2} =  \zeta_3 ( 9 c_j +2 \sqrt{2} d_j + e_j ), &
		b_{{j+8},3} & = \mathrm{i} (8\sqrt{2} c_j + 4 d_j + \sqrt{2} e_j + h_j ),
	\end{aligned}$$
	%\begin{align*}
	%	a_{{j+8},\,4} & =\zeta_3 a_{{j+8},\,0} = \zeta_3\,  a_j,\\
	%	a_{{j+8},\,3} & =\zeta_6 a_{{j+8},\,1} = \zeta_6\, 2 \sqrt{2}\, a_j+b_j,\\
	%	a_{{j+8},\,2} & = \zeta_6\, ( a_j + \sqrt{2}\, b_j + g_j),\\
	%	b_{{j+8},\,6} & = -b_{{j+8},\,0} = c_j,\\
	%	b_{{j+8},\,5} & =b_{{j+8},\,1}=\zeta_{12}^5 (3 \sqrt{2} c_j + d_j),\\
	%	b_{{j+8},\,4} & =b_{{j+8},\,2} =  \zeta_3 ( 9 c_j +2 \sqrt{2} d_j + e_j ),\\
	%	b_{{j+8},\,3} & = \mathrm{i} (8\sqrt{2} c_j + 4 d_j + \sqrt{2} e_j + h_j ),
	%\end{align*}
	for $j=1,\ldots, 8$. 
	
	We note that the points $P'_j= \left(  \zeta_{12}^2 t^2 x(P_j),  \zeta_{12}^3 t^3  y(P_j) \right)$ belong to
	the Mordell--Weil lattice of $\Ee: y^2= x^3+ t^{12}+1 $, which is birational to $\Ee_6$	over $\CC(t)$.
	See \cite{Usui2008} for more details.	
	Having  polynomial coordinates, the  $P'_j$'s have no intersection with zero sections of $\Ee$, we get that
	$\langle P'_j, P'_j \rangle = 4$,
	$\langle P'_{j_1}, P'_{j_2}\rangle = 2- (P'_{j_1} \cdot P'_{j_2})$,
	and  for any $1\leq j_1 \neq j_2 \leq 16$, the intersection number
	$(P'_{j_1} \cdot P'_{j_2})$ can be computed by:
	\begin{align*}
		( P'_{j_1} \cdot P'_{j_2}) &= \deg(\gcd (x_{j_1}- x_{j_2} , y_{j_1}- y_{j_2})) \nonumber \\ 
		&\quad + \min \{4- \deg( x_{j_1}- x_{j_2}) , 6- \deg( y_{j_1}- y_{j_2})\}.
	\end{align*}
	Thus,  we obtain the Gram matrix $R_6$ with determinant is  $2^4 \cdot 3^4$ of the height pairing  for  $P'_j$'s and hence $P_j$'s:
	\[
	\tiny % Reduces font size to fit page width
	\setcounter{MaxMatrixCols}{20} % Allows 16 columns
	R_6=
	\begin{pmatrix}
		4 & 2 & 0 & 0 & 0 & 2 & -1 & 2 & 0 & 0 & 0 & 0 & 0 & 0 & 0 & 0 \\ \noalign{\medskip}
		2 & 4 & 1 & 1 & 1 & 1 & -2 & 1 & 0 & 0 & 0 & 0 & 0 & 0 & 0 & 0 \\ \noalign{\medskip}
		0 & 1 & 4 & 0 & 0 & -2 & 0 & 0 & 0 & 0 & -1 & 1 & 0 & 0 & 0 & 0 \\ \noalign{\medskip}
		0 & 1 & 0 & 4 & -2 & 0 & -2 & 1 & 1 & 0 & 0 & 0 & 0 & 0 & 0 & 0 \\ \noalign{\medskip}
		0 & 1 & 0 & -2 & 4 & 0 & 1 & -2 & -1 & 0 & 0 & 0 & 0 & 0 & 0 & 0 \\ \noalign{\medskip}
		2 & 1 & -2 & 0 & 0 & 4 & 0 & 0 & 0 & 0 & 0 & -1 & 0 & 0 & 0 & 0 \\ \noalign{\medskip}
		-1 & -2 & 0 & -2 & 1 & 0 & 4 & -2 & -1 & 1 & 0 & 0 & 0 & 0 & 0 & 0 \\ \noalign{\medskip}
		2 & 1 & 0 & 1 & -2 & 0 & -2 & 4 & 0 & -1 & 0 & 0 & 0 & 0 & 0 & 0 \\ \noalign{\medskip}
		0 & 0 & 0 & 1 & -1 & 0 & -1 & 0 & 4 & -2 & 0 & 0 & 0 & 2 & 0 & 2 \\ \noalign{\medskip}
		0 & 0 & 0 & 0 & 1 & 0 & 1 & -1 & -2 & 4 & 0 & 0 & -2 & 0 & -2 & 0 \\ \noalign{\medskip}
		0 & 0 & -1 & 0 & 0 & 0 & 0 & 0 & 0 & 0 & 4 & -2 & 0 & 2 & 0 & -2 \\ \noalign{\medskip}
		0 & 0 & 1 & 0 & 0 & -1 & 0 & 0 & 0 & 0 & -2 & 4 & -2 & 0 & 2 & 0 \\ \noalign{\medskip}
		0 & 0 & 0 & 0 & 0 & 0 & 0 & 0 & 0 & -2 & 0 & -2 & 4 & -2 & 0 & 0 \\ \noalign{\medskip}
		0 & 0 & 0 & 0 & 0 & 0 & 0 & 0 & 2 & 0 & 2 & 0 & -2 & 4 & 0 & 0 \\ \noalign{\medskip}
		0 & 0 & 0 & 0 & 0 & 0 & 0 & 0 & 0 & -2 & 0 & 2 & 0 & 0 & 4 & -2 \\ \noalign{\medskip}
		0 & 0 & 0 & 0 & 0 & 0 & 0 & 0 & 2 & 0 & -2 & 0 & 0 & 0 & -2 & 4
	\end{pmatrix}
	\]		
	Thus, the points $P_j$'s, $j=1, \dots, 16$, form a set of independent generators of $\mathcal{E}_6$ over $\mathcal{K}_6(t)$.
	We refer the reader to see  \cite[check-6]{k3-codes}  for all computations in this section.
	\iffalse
	$$	\left[\begin{array}{cccccccccccccccc}
		4 & 2 & 0 & 0 & 2 & 0 & 0 & 0 & 0 & 0 & 0 & 0 & 0 & 0 & 0 & 0 
		\\
		2 & 4 & 0 & 0 & 2 & 0 & 0 & 0 & 0 & 0 & 0 & 0 & 0 & 0 & 0 & 0 
		\\
		0 & 0 & 4 & 0 & 0 & -2 & 0 & 2 & 0 & 0 & 0 & 0 & 0 & -1 & 0 & 0 
		\\
		0 & 0 & 0 & 4 & 0 & 2 & 2 & 0 & 0 & 0 & 0 & 0 & 0 & 1 & 1 & 1 
		\\
		2 & 2 & 0 & 0 & 4 & 0 & 0 & 2 & 0 & 0 & 0 & 0 & 0 & 0 & 2 & 1 
		\\
		0 & 0 & -2 & 2 & 0 & 4 & 2 & 0 & 0 & 0 & 1 & -1 & 0 & 0 & 0 & 1 
		\\
		0 & 0 & 0 & 2 & 0 & 2 & 4 & 2 & 0 & 0 & 0 & -1 & -2 & 0 & 0 & -1 
		\\
		0 & 0 & 2 & 0 & 2 & 0 & 2 & 4 & 0 & 0 & 0 & -1 & -1 & -1 & 1 & 0 
		\\
		0 & 0 & 0 & 0 & 0 & 0 & 0 & 0 & 4 & 2 & 0 & 0 & 2 & 0 & 0 & 0 
		\\
		0 & 0 & 0 & 0 & 0 & 0 & 0 & 0 & 2 & 4 & 0 & 0 & 2 & 0 & 0 & 0 
		\\
		0 & 0 & 0 & 0 & 0 & 1 & 0 & 0 & 0 & 0 & 4 & 0 & 0 & -2 & 0 & 2 
		\\
		0 & 0 & 0 & 0 & 0 & -1 & -1 & -1 & 0 & 0 & 0 & 4 & 0 & 2 & 2 & 0 
		\\
		0 & 0 & 0 & 0 & 0 & 0 & -2 & -1 & 2 & 2 & 0 & 0 & 4 & 0 & 0 & 2 
		\\
		0 & 0 & -1 & 1 & 0 & 0 & 0 & -1 & 0 & 0 & -2 & 2 & 0 & 4 & 2 & 0 
		\\
		0 & 0 & 0 & 1 & 2 & 0 & 0 & 1 & 0 & 0 & 0 & 2 & 0 & 2 & 4 & 2 
		\\
		0 & 0 & 0 & 1 & 1 & 1 & -1 & 0 & 0 & 0 & 2 & 0 & 2 & 0 & 2 & 4 
	\end{array}\right]
	$$
	\fi
	\bigskip 

	\bibliographystyle{amsplain}
	
	\bibliography{SBIB1}{}

\end{document}